\documentclass[aop,MSNbibl,seceqn,dvips]{arximspdf}

\doi{10.1214/09-AOP510}
\volume{38}
\issue{3}
\pubyear{2010}
\firstpage{1320}
\lastpage{1344}

\makeatletter

\newtheorem{theorem}{Theorem}[section]
\newtheorem{lemma}{Lemma}[section]

\def\eqref#1{(\ref{#1})}
\makeatother

\begin{document}
\begin{frontmatter}

\title{Large faces in Poisson hyperplane mosaics}
\runtitle{Poisson hyperplane mosaics}

\begin{aug}
\author[A]{\fnms{Daniel} \snm{Hug}\corref{}\ead[label=e1]{daniel.hug@kit.edu}} and
\author[B]{\fnms{Rolf} \snm{Schneider}\ead[label=e2]{rolf.schneider@math.uni-freiburg.de}}

\runauthor{D. Hug and R. Schneider}

\affiliation{Karlsruhe Institute of Technology and
Albert--Ludwigs-Universit\"{a}t Freiburg}

\address[A]{Department of Mathematics\\Karlsruhe Institute of
Technology\\
D-76128 Karlsruhe\\Germany\\ \printead{e1}}
\address[B]{Mathematisches Institut\\Albert--Ludwigs-Universit\"{a}t Freiburg\\
Eckerstr. 1\\D-79104 Freiburg i. Br.\\Germany\\ \printead{e2}\\}

\end{aug}

\received{\smonth{5} \syear{2009}}
\revised{\smonth{10} \syear{2009}}

%
\begin{abstract}
A generalized version of a well-known problem of D. G. Kendall states
that the zero cell of a stationary Poisson hyperplane
tessellation in ${\mathbb R}^d$, under the condition that it has large
volume, approximates with high probability a certain
definite shape, which is determined by the directional distribution of
the underlying hyperplane process.
This result is extended here to typical $k$-faces of the tessellation,
for $k\in\{2,\dots,d-1\}$.
This requires the additional condition that the direction of the face
be in a sufficiently small neighbourhood of a given direction.
\end{abstract}

%


\end{frontmatter}

\section{Introduction}

A well-known problem of D. G. Kendall, popularized in the foreword to the
first edition (1987) of the book \cite{SKM95}, asked whether the shape
of the zero cell of a stationary, isotropic Poisson line process in the
plane, under the condition that the cell has large area, must be
approximately circular, with high probability. An affirmative answer
was given by Kovalenko \cite{Kov97,Kov99}. Several higher-dimensional
versions and variants of Kendall's problem were treated in \cite{HRS04a,HRS04b,HS04,HS05,HS07a,HS07b}. In \cite{HRS04a}, the
subject of investigation was the zero cell of a stationary Poisson
hyperplane process with a general (nondegenerate) directional
distribution in $d$-dimensional Euclidean space, under the condition
that the cell has large volume. The asymptotic shape of such cells was
found to be that of the so-called Blaschke body of the hyperplane
process. This is (up to a dilatation) the convex body, centrally
symmetric with respect to the origin, that has the spherical
directional distribution of the hyperplane process as its surface area
measure. Its existence and uniqueness follow from a celebrated theorem
going back to Minkowski.

The purpose of the present paper is an extension of the latter result
to $k$-dimensional faces, for $k\in\{2,\dots,d-1\}$. The natural
extension of the zero cell, which is stochastically equivalent to the
volume weighted typical cell, is the notion of the
($k$-volume-)weighted typical $k$-face. We consider the weighted typical $k$-face
under the condition that it has large $k$-dimensional volume and that
its direction space (the translate of its affine hull passing through
the origin) is in a small neighbourhood of a given $k$-dimensional
subspace $L^*$. We can then again identify an asymptotic shape, namely
that of the Blaschke body of the section process of the given
hyperplane process with the subspace $L^*$. The main results, whose
precise formulation requires some preparations, are formulated in the
theorems at the end of the next section. The extension from cells to
lower-dimensional faces is not routine; the proof has become possible
through a recently established representation for the distribution of
the weighted typical $k$-face (\cite{Sch09}, Theorem 1) and a special
stability result for the convex bodies obtained from Minkowski's
existence theorem, which was proved in \cite{HS02}.

Once the result is proved for weighted typical $k$-faces (Theorem \ref
{T1}), it can be used to derive a variant for typical $k$-faces
(Theorem \ref{T2}). From these theorems, the existence of limit shapes
can be deduced (Theorem~\ref{T7}).

\section{Preliminaries and main results}\label{sec2}

Fundamental facts about Poisson hyperplane processes and random
mosaics, as well as corresponding notions that are not explained here,
can be found in the book \cite{SW08}. For the employed notions and
results from convex geometry, we refer to \cite{Sch93}.

We denote by $\mathbb{R}^d$ the $d$-dimensional Euclidean vector space
(assuming $d\ge3$ throughout), with scalar product $\langle\cdot
,\cdot\rangle$ and induced norm $\|\cdot\|$. Its unit ball and unit
sphere are denoted by ${\mathbb B}^d$ and ${\mathbb S}^{d-1}$,
respectively. Further, $\mathit{SO}_d$ is the rotation group, $G(d,k)$ is the
Grassmannian of $k$-dimensional linear subspaces of $\mathbb{R}^d$ and
$A(d,k)$ is the set of $k$-flats ($k$-dimensional affine subspaces) of
$\mathbb{R}^d$; all these sets are equipped with their standard topologies.\vspace*{1pt}

By $\mathcal{H}^d\subset A(d,d-1)$, we denote the space of hyperplanes
in $\mathbb{R}
^d$ not passing through the origin $\mathbf{o}$. Every hyperplane in
$\mathcal{H}^d$
has a unique representation
\[
H(\mathbf{u},t)=\{\mathbf{x}\in\mathbb{R}^d\dvtx \langle\mathbf
{x},\mathbf{u}\rangle=t\}
\]
with $\mathbf{u}\in{\mathbb S}^{d-1}$ and $t>0$, and
\[
H^-(\mathbf{u},t)=\{\mathbf{x}\in\mathbb{R}^d\dvtx \langle\mathbf
{x},\mathbf{u}\rangle\le t\}
\]
is the closed halfspace bounded by it that contains $\mathbf{o}$. We write
$H^-= H^-(\mathbf{u},t)$ if $H=H(\mathbf{u},t)$.

Let $\mathcal{K}$ be the space of convex bodies (nonempty, compact, convex
subsets) in ${\mathbb R}^d$, endowed with the Hausdorff metric $\delta$.
For $k\in\{2,\dots,d-1\}$ and a subspace $L\in G(d,k)$, we denote by
$\mathcal{K}(L)$ the set of convex bodies $K\subset L$ and by $\mathcal{
K}_0(L)$ the subset of $k$-dimensional bodies $K$ with $\mathbf{o}\in
{\rm
relint}  K$ (where relint denotes the relative interior).

In the following, measures on a given topological space $T$, if not
further specified, are always positive measures on the Borel $\sigma
$-algebra $\mathcal{B}(T)$ of the space.

We turn to hyperplane processes. As usual and convenient in the theory
of point processes, we often identify a simple counting measure $\eta$
on a topological space $E$ with its support, so that $\eta(\{x\})=1$
and $x\in\eta$ are used synonymously, and $\eta(A)$ and $\operatorname{card}(\eta\cap A)$ both denote the number of elements of $\eta$ in
the subset $A\subset E$.

Let $X$ be a stationary Poisson hyperplane process in $\mathbb{R}^d$. We
denote the underlying probability by $\mathbf{P}$ and mathematical
expectation by $\mathbf{E}$. The intensity measure $\Theta= \mathbf
{E}X(\cdot)$
of $X$ has a representation (equivalent to \cite{SW08}, (4.33))
\[
\Theta(A) = 2\gamma\int_{{\mathbb S}^{d-1}}\int_0^\infty\mathbf{1}
_A(H(\mathbf{u},t)) \,\mathrm{d}t \,\varphi(\mathrm{d}\mathbf{u})
\]
for $A\in\mathcal{B}(A(d,d-1))$, where $\gamma$ is the intensity of $X$
and $\varphi$ is its spherical directional distribution. This is an
even probability measure on the unit sphere; we assume that it is not
concentrated on any great subsphere.

Together with the hyperplane process $X$, the following processes of
lower-dimensional flats derived from it will play an essential role.
First, let $k\in\{2,\dots,d-1\}$ and $L\in G(d,k)$. The section
process $X\cap L$ is obtained by taking all $(k-1)$-dimensional
intersections of hyperplanes
of $X$ with $L$; see \cite{SW08}, pages~129 ff. It is a stationary
Poisson process of $(k-1)$-flats in $L$. We denote its intensity by
$\gamma_{X\cap L}$ and its spherical directional distribution, defined
on ${\mathbb S}^{d-1}\cap L$, by $\varphi_{X\cap L}$. Second, for
$k\in\{0,\dots,d-1\}$, the process $X_{d-k}$ is obtained by
intersecting any $d-k$ hyperplanes of $X$ which are in general
position; see \cite{SW08}, Section~4.4. It~is a stationary process of
$k$-flats and is called the intersection process of order $d-k$ of $X$.
We denote its intensity by $\gamma_{d-k}$ and its directional
distribution by $\mathbf{ Q}_{d-k}$. The latter is a probability measure on
$G(d,k)$.

The hyperplane process $X$ induces a tessellation $X^{(d)}$ of $\mathbb{R}^d$
and with it the process $X^{(k)}$ of its $k$-dimensional faces, for
$k=0,\dots,d-1$ (for the notation, note the slight digression from
\cite{SW08}, where $X$ and $X^{(d)}$ are denoted by $\widehat X$ and
$X$, resp.). The zero cell of $X^{(d)}$ is the cell ($d$-face)
containing $\mathbf{o}$ and thus is the random polytope given by
\[
Z_0:= \bigcap_{H\in X} H^-.
\]
Its counterpart for $k$-faces can be defined as follows (see \cite
{BL07,Sch09}, e.g.). Let $M_k$ denote the random
measure defined by restricting the $k$-dimensional Hausdorff measure to
the union of the $k$-flats of $X_{d-k}$. Further, let ${\sf N}_s$
denote the set of simple counting measures on $A(d,d-1)$
and $\mathcal{N}_s$ the usual $\sigma$-algebra of ${\sf N}_s$ (see
\cite{SW08}, Section 3.1).
Let $B\subset\mathbb{R}^d$ be a Borel set of Lebesgue measure $1$ and let
$\mathcal{A}\in\mathcal{N}_s$. Then
\[
\mathbf{P}^0_k(\mathcal{A}):= \frac{1}{\mathbf{E}M_k(B)} \mathbf
{E}\int_B \mathbf{1}
_{\mathcal{A}}(X-\mathbf{x}) M_k(\mathrm{d}\mathbf{x})
\]
defines a probability measure $\mathbf{P}^0_k$ (a Palm distribution,
independent of $B$) on the measurable space $({\sf N}_s,\mathcal{N}_s)$.
Let $Y$ be a hyperplane process with distribution $\mathbf{P}^0_k$.
Then the
\textit{weighted typical $k$-face} $Z_0^{(k)}$ of $X$ (i.e., of the
mosaic induced by $X$) is defined as the a.s. unique $k$-face in
$Y^{(k)}$ containing the origin $\mathbf{o}$. The distribution of the random
polytope $Z_0^{(k)}$ is uniquely determined and coincides, up to
translations, with that of the typical $k$-face $Z^{(k)}$ weighted by
its $k$-volume. This is revealed by the relation
%
\begin{equation}\label{1a}
\mathbf{E}f\bigl(Z_0^{(k)}\bigr) = \frac{1}{\mathbf{E}V_k(Z^{(k)})} \mathbf
{E}
\bigl[f\bigl(Z^{(k)}\bigr)V_k\bigl(Z^{(k)}\bigr)\bigr],
\end{equation}
holding for every translation invariant, nonnegative, measurable
function $f$ on the space of $k$-dimensional polytopes (see
\cite{Sch09}, equation (11)). Here, $V_k$ denotes the $k$-dimensional volume, and
$Z^{(k)}$ is the typical $k$-face of the mosaic induced by $X$, as
defined in \cite{SW08}, page 450. An even more intuitive interpretation
of the weighted typical $k$-face, up to translations, is the following.
Let $s$ denote the Steiner point (or any other centre function, see
\cite{SW08}, page 110), and let $W\in\mathcal{K}$ be an arbitrary convex
body with positive volume. Then, for every Borel set $A$ in the space
of convex polytopes,
\[
\mathbf{P}\bigl\{Z_0^{(k)}-s\bigl(Z_0^{(k)}\bigr)\in A\bigr\} = \lim_{r\to
\infty}
\frac{\mathbf{E}\sum_{F\in X^{(k)}, F\subset rW} \mathbf
{1}_A(F-s(F))V_k(F)} {\mathbf{E}
\sum_{F \in X^{(k)}, F\subset rW} V_k(F)}.
\]

The following integral representation for the distribution of
$Z_0^{(k)}$ is proved in \cite{Sch09}, Theorem 1. For Borel sets $A$
in the space of convex polytopes,
%
\begin{equation}\label{n1}
\mathbf{P}\bigl\{Z_0^{(k)}\in A\bigr\} = \int_{G(d,k)} \mathbf
{P}\{Z_0\cap L\in
A\} \mathbf{ Q}_{d-k}(\mathrm{d}L).
\end{equation}

Recall (\cite{SW08}, page 162) that the \textit{Blaschke body} of $ X$ is
the $\mathbf{o}$-symmetric convex body $B( X)$ with surface area measure
$S_{d-1}(B( X),\cdot) = \gamma\varphi$. To describe the asymptotic
shape of large weighted typical $k$-faces, we need the Blaschke body
$B(X\cap L)$ of the section process $X\cap L$, for $L$ in the support
of the measure $\mathbf{ Q}_{d-k}$. Since only the homothety class of the
Blaschke body plays a role in the following, we may replace it by any
dilate. It is convenient here to use the $\mathbf{o}$-symmetric body
$B_L\subset L$ with surface area measure on ${\mathbb S}^{d-1}\cap L$
given by the spherical directional distribution $\varphi_{X\cap L}$ of
the section process $X\cap L$.

We need some particular notions of distance.

For a rotation $\rho\in \mathit{SO}_d$, let $M_\rho$ be the matrix of $\rho$
with respect to the standard orthonormal basis of ${\mathbb R}^d$. We
define the distance of $\rho$ from the identity by
\[
|\rho| := \|M_\rho-I\|,
\]
where $I$ is the unit matrix and $\|A\|=(\sum_{i,j=1}^d
a_{ij}^2)^{1/2}$ is the Frobenius norm of the matrix $A=
(a_{ij})_{i,j=1}^d$. Note that $|\rho|=|\rho^{-1}|$, since $M_{\rho
^{-1}}-I$ is the transpose of $M_{\rho}-I$, and that for $\mathbf
{x}\in
{\mathbb R}^d$ we have
\[
\|\mathbf{x}-\rho\mathbf{x}\| \le|\rho|\|\mathbf{x}\|.
\]

On $G(d,k)$, we introduce a metric $\Delta$ by
\[
\Delta(L,E) := \min\{|\rho|\dvtx \rho\in \mathit{SO}_d,  \rho L=E\}.
\]
The triangle inequality follows from
\[
\|M_{\rho_1}M_{\rho_2}-I\| \le\|M_{\rho_1}M_{\rho_2}-M_{\rho_1}\| +
\|M_{\rho_1}-I\| = \|M_{\rho_2}-I\| + \|M_{\rho_1}-I\|
\]
for $\rho_1,\rho_2\in \mathit{SO}_d$. The metric $\Delta$ induces the
standard topology of $G(d,k)$ and is particularly convenient for us.
For metrics on Grassmannians involving, like this one, a ``direct
rotation'' between subspaces, we refer to \cite{DK70}, the survey
article \cite{PW94}, and the references given there.

For $\theta>0$, the $\theta$-neighbourhood of a subspace $L^*\in
G(d,k)$ is defined by
\[
N_\theta(L^*):= \{L\in G(d,k)\dvtx  \Delta(L,L^*)<\theta\}.
\]

For $L\in G(d,k)$ and $K,M\in\mathcal{K}_0(L)$ with $M=-M$, let
%
\begin{equation}\label{n2}
\vartheta(K,M) :=\log\min\{\beta/\alpha\dvtx  \alpha,\beta>0,
\exists \mathbf{z}\in L\dvtx \alpha M \subset K+\mathbf{z}\subset\beta
M\}.
\end{equation}
The function $\vartheta$ measures the deviation of the homothetic
shapes of $K$ and $M$; it is nonnegative, and it vanishes if and only
if $K$ and $M$ are homothetic.

For $L,E\in G(d,k)$ and convex bodies $K\in\mathcal{K}_0(L)$ and $M\in
\mathcal{K}_0(E)$ with $M=-M$, let
%
\begin{equation}\label{n3}
\vartheta(K,M) := \min\{\vartheta(\rho K,M)\dvtx \rho\in \mathit{SO}_d, \rho
L=E, |\rho| = \Delta(L,E)\}.
\end{equation}
Note that this definition is consistent with (\ref{n2}), since $|\rho
|= \Delta(L,E)$ in the case $L=E$ implies that $\rho$ is the
identity. Note also that $\vartheta(K,M)$ is symmetric in $K$ and $M$
if both bodies are $\mathbf{o}$-symmetric.

For a $k$-dimensional convex body $K$, we denote by $D(K)=\operatorname{lin}(K-K)\in G(d,k)$ its direction space; this is the linear subspace
parallel to the affine hull of $K$.

Throughout the paper, several constants $c_i$ will appear, which may
depend on various data. Their possible dependence on the dimension $d$
will not be mentioned, since we work in a space of fixed dimension.

Now, we can formulate our main result.

\begin{theorem}\label{T1} Let $X$ be a stationary Poisson hyperplane
process in $\mathbb{R}^d$ with intensity $\gamma$ and spherical directional
distribution $\varphi$. Let $k \in\{2,\dots,d-1\}$, and let
$Z_0^{(k)}$ be the weighted typical $k$-face of the mosaic induced by
$X$. Let $\mathbf{ Q}_{d-k}$ be the directional distribution of the
intersection process of order $d-k$ of $X$.

Let $\varepsilon>0$ be given. Then there exist constants $c_1,c_2>0$,
depending only on $\varphi, \gamma, \varepsilon$, and a constant
$c_3>0$, depending only on $\varphi, \gamma$, such that the following
is true.
If $L^* \in G(d,k)$ is in the support of the measure $\mathbf{ Q}_{d-k}$, then
\begin{eqnarray*}
 &&\mathbf{P}\bigl\{ \vartheta\bigl(Z_0^{(k)},B_{L^*}\bigr) \ge\varepsilon\vert
V_k\bigl(Z_0^{(k)}\bigr) \ge a,  D\bigl(Z_0^{(k)}\bigr) \in N_{\theta}(L^*)\bigr\}
\\
&&\qquad \le c_2 \exp[-c_3 \varepsilon^{k+1}a^{1/k}]
\end{eqnarray*}
for all $a\ge1$ and all $0<\theta\le c_1$.
\end{theorem}

In other words, if a subspace $L^*$ in the support of the distribution
$\mathbf{ Q}_{d-k}$ and a bound $\varepsilon>0$ are given, then the
probability that the weighted typical cell $Z_0^{(k)}$ deviates in
shape from the (dilated) Blaschke body $B_{L^*}$ by at least $\varepsilon
$, under the condition that its direction space is contained in a
suitable neighbourhood of $L^*$ and its volume is at least $a>0$,
becomes exponentially small for large $a$. From this, one can deduce
that the Blaschke body is the limit shape of $Z_0^{(k)}$ if the volume
of $Z_0^{(k)}$ tends to infinity and its direction space tends to $L^*$
(see Theorem \ref{T7} for a precise formulation). The assumptions of
Theorem \ref{T1} are inevitable: the subspace $L^*$ must be chosen in
the support of the measure $\mathbf{ Q}_{d-k}$ since, by (\ref{n1}), the
direction space of $Z_0^{(k)}$ lies almost surely in the support of
this measure. (This is also intuitively obvious: the $k$-faces of the
tessellation $X^{(d)}$ are generated by intersections of hyperplanes
from the process $X$.) Further, the Blaschke body $B_{L^*}$ depends on
$L^*$, hence in general only weighted typical cells with a direction
space close to $L^*$ can approximate the shape of $B_{L^*}$. The
admissible size of the neighbourhood $N_\theta(L^*)$ in Theorem \ref
{T1} depends heavily on Lemma~\ref{L4} below and thus on the
underlying stability theorem for Minkowski's existence theorem.
Lemma~\ref{L4} would yield additional information on the dependence of $c_1$ on
$\varepsilon$, but a more explicit specification of the neighbourhood
will only be possible for directional distributions $\varphi$ where
the solutions of Minkowski's problem are more explicitly accessible.

There are, however, two simple cases which should be mentioned. If the
hyperplane process $X$ is isotropic, that is, its directional
distribution $\varphi$ is invariant under rotations, then all Blaschke
bodies $B_{L^*}$ are balls, and the condition on the direction space
$D(Z_0^{(k)})$ can be omitted entirely. In fact, if $X$ is isotropic,
then it can be deduced from (\ref{n1}) (see \cite{Sch09}) that there
exists a random rotation $\rho$ such that $\rho Z_0^{(k)}$ has the
same distribution as the zero cell of a stationary isotropic Poisson
$(k-1)$-flat process in a fixed $k$-dimensional subspace of ${\mathbb
R}^d$. Therefore, one can immediately apply the results from \cite
{HRS04a} in that subspace.

Another simple case is that of a discrete directional distribution. If
the directional distribution $\varphi$ of $X$ is concentrated in
finitely many points, then every body $B_{L^*}$ is a $k$-dimensional
polytope. The distribution $\mathbf{ Q}_{d-k}$ is concentrated in finitely
many elements of $G(d,k)$. Hence, there exist only finitely many
possibilities for the direction space of $Z_0^{(k)}$. If $L^*$ in the
support of $\mathbf{ Q}_{d-k}$ is given, one can then choose for $N_\theta
(L^*)$ in Theorem \ref{T1} any neighbourhood of $L^*$ containing no other
element of the support of $\mathbf{ Q}_{d-k}$.

The proof of Theorem \ref{T1} will be given in Section \ref{sec5}. The next
section provides geometric results in preparation for that proof.

The arguments leading to Theorem \ref{T1} can be modified to yield
also a corresponding result for the typical $k$-face $Z^{(k)}$ of the
mosaic induced by $X$.

\begin{theorem}\label{T2} The assertion of Theorem \ref{T1}
remains true if the weighted typical $k$-face
$Z_0^{(k)}$ is replaced by the typical $k$-face $Z^{(k)}$ of the mosaic
induced by~$X$.
\end{theorem}

This is in analogy to the corresponding result for the typical cell,
Theorem 2 in~\cite{HRS04a}.

We have restricted ourselves here, in agreement with D. G. Kendall's original
question, to the volume functional. For the zero cell, we have
investigated in \cite{HS07a} asymptotic shapes when the size of the
zero cell is measured by various other functionals. It is a natural
question whether such results carry over to $k$-faces. This is
certainly possible in the isotropic case and for rotation invariant
size functionals, by the remark made above. However, for non-isotropic
distributions and general size functionals, asymptotic shapes are no
longer controlled by the Blaschke body, so that the crucial Lemma \ref{L4}
below must be replaced by a different approach.

\section{Auxiliary continuity and stability results}\label{sec3}

Throughout this paper, $ X$ is a stationary Poisson hyperplane process
in ${\mathbb R}^d$, with intensity $\gamma$ and spherical directional
distribution $\varphi$. We assume that $\varphi$ is not concentrated
on a great subsphere and, without loss of generality, that it is even
(invariant under reflection in $\mathbf{o}$). The main topic of this
section is
the dependence of the dilated Blaschke bodies $B_L\subset L$, $L\in
G(d,k)$, on the probability measure $\varphi$ and on $L$.

Let $L\in G(d,k)$, where $k\in\{2,\dots,d-1\}$. The set ${\mathbb
S}_L^{k-1}:= {\mathbb S}^{d-1}\cap L$ is the unit sphere in $L$. The
surface area measure of $K\in\mathcal{K}_0(L)$ is denoted by
$S^L_{k-1}(K,\cdot)$; this is a measure on ${\mathbb S}^{k-1}_L$. By
definition, the Blaschke body $B(X\cap L)$ is the unique convex body in
$\mathcal{K}(L)$, centrally symmetric with respect to $\mathbf{o}$, for which
\[
S^L_{k-1}\bigl(B(X\cap L),\cdot\bigr) =\gamma_{X\cap L}\varphi_{X\cap L},
\]
where $\gamma_{X\cap L}$ is the intensity and $\varphi_{X\cap L}$ is
the spherical directional distribution of the section process $X\cap
L$. Existence and uniqueness of this body follow from Minkowski's
theorem (see \cite{Sch93}, Section 7.1, e.g.). We work here with
a dilate of the Blaschke body, the $\mathbf{o}$-symmetric body $B_L$
defined by
\[
S^L_{k-1}(B_L,\cdot) =\varphi_{X\cap L}.
\]

The \textit{associated zonoid} $\Pi_{ X}$ of $ X$ is the projection body
of $B( X)$ and hence has generating measure $\frac{1}{2}\gamma\varphi
$, that is, its support function has the integral representation
\[
h(\Pi_X, \mathbf{u}) = V_{d-1}(B(X)|\mathbf{u}^\perp)=\frac{\gamma
}{2} \int
_{{\mathbb S}^{d-1}} |\langle\mathbf{u},\mathbf{v}\rangle| \varphi
(\mathrm{d}\mathbf{v}
),\qquad\mathbf{u}\in{\mathbb S}^{d-1},
\]
where $\cdot|L$ denotes the orthogonal projection to $L$ and $\mathbf{u}
^\perp$ is the hyperplane through $\mathbf{o}$ orthogonal to $\mathbf
{u}$. In the
following, the support function $h(K,\cdot)$ of a convex body $K$ is
always defined on $\mathbb{R}^d$, also if $K\subset L$, $L\in G(d,k)$.

Let $L\in G(d,k)$. The associated zonoid of the section process $ X\cap
L$ is given by
\[
\Pi_{ X\cap L}= \Pi_{ X}|L
\]
(see \cite{SW08}, equation (4.61)). From this, we can read off the
generating measure of the zonoid $\Pi_{ X\cap L}$ and hence the
essential parameters of the section process $ X\cap L$. We define the
spherical projection ${\rm pr}_L\dvtx {\mathbb S}^{d-1}\setminus L^\perp\to
{\mathbb S}_L^{k-1}$ by
\[
\operatorname{pr}_L(\mathbf{u}) := \frac{\mathbf{u}|L}{\|\mathbf{u}|L\|}
\qquad\mbox{for } \mathbf{u}
\in{\mathbb S}^{d-1}\setminus L^\perp.
\]
Let $\mathcal{M}({\mathbb S}^{d-1})$ denote the cone of finite Borel
measures on ${\mathbb S}^{d-1}$. The spherical projection $\pi_L\dvtx \mathcal{
M}({\mathbb S}^{d-1})\to\mathcal{M}({\mathbb S}_L^{k-1})$ is defined by
%
\begin{equation}\label{0}
\pi_L\mu(A) = \int_{{\mathbb S}^{d-1}\setminus L^\perp} \mathbf{
1}_A(\operatorname{pr}_L(\mathbf{u}))\|\mathbf{u}|L\| \mu(\mathrm{d}\mathbf{u})
\end{equation}
for Borel sets $A\subset{\mathbb S}_L^{k-1}$ and for $\mu\in\mathcal{
M}({\mathbb S}^{d-1})$. (More general spherical projections and their
applications are treated in \cite{GKW08}.)

For a segment $S= \operatorname{conv}\{-\alpha\mathbf{v},\alpha\mathbf{v}\}$
with $\mathbf{v}\in
{\mathbb S}^{d-1}\setminus L^\perp$ and $\alpha>0$, we have for
$\mathbf{u}
\in\mathbb{R}^d$,
\[
h(S|L,\mathbf{u}) = \bigl|\langle\mathbf{v}|L,\mathbf{u}\rangle\bigr|\alpha
= |\langle{\rm
pr}_L(\mathbf{v}),\mathbf{u}\rangle|\cdot\|\mathbf{v}|L\|\alpha.
\]
Hence, if $Z$ is a zonoid with generating measure $\mu$, then the
zonoid $Z|L$ has generating measure $\pi_L\mu$. In particular, the
generating measure of $\Pi_{ X}|L$ is given by $\frac{1}{2}\gamma\pi
_L\varphi$. It follows that the Blaschke body $B(X\cap L)$ has surface
area measure $S^L_{k-1}(B( X\cap L),\cdot) = \gamma\pi_L\varphi$.
We conclude that
%
\begin{equation}\label{5a}
\gamma_{X\cap L}\varphi_{X\cap L}= \gamma\pi_L\varphi.
\end{equation}
This is \cite{SW08}, Theorem 4.4.7, for hyperplane processes and in
terms of spherical directional distributions.

Since the bodies $B_L$ are obtained from the (nonconstructive)
existence theorem of Minkowski, it is not trivial that they depend
continuously on $L$. We need a stronger result, estimating how close
$B_L$ and $B_E$ are in a suitable sense if the subspaces $L,E\in
G(d,k)$ are close to each other. Such an estimate (Lemma \ref{L4}) is
obtained from a stability result for Minkowski's theorem that uses the
Prokhorov metric for measures. (Diskant's stability result (see \cite
{Sch93}, Theorem 7.2.), which is in terms of the total variation norm
of the difference, would not be strong enough for this purpose.) For
finite measures
$\mu,\nu$ on ${\mathbb S}^{d-1}$, the Prokhorov distance $d_P(\mu
,\nu)$ is defined by
\begin{eqnarray*}
d_P(\mu,\nu) &=& \inf\{\varepsilon>0\dvtx \mu(A)\le\nu(A_\varepsilon
)+\varepsilon\mbox{ and } \nu(A)\le\mu(A_\varepsilon)+\varepsilon
\\
&&{}\hspace*{96pt}\quad\ \ \mbox{for all Borel sets }A\subset{\mathbb S}^{d-1}\},
\end{eqnarray*}
where
\[
A_\varepsilon:= \{ \mathbf{y}\in{\mathbb S}^{d-1}\dvtx \|\mathbf{x}-\mathbf
{y}\|<\varepsilon\mbox{
for some }\mathbf{x}\in A\}.
\]
Analogous definitions are used for measures on ${\mathbb S}^{k-1}_L$.
For a rotation $\rho$ and a measure $\mu$, we denote by $\rho\mu$
the image measure of $\mu$ under $\rho$, defined by $(\rho\mu
)(A)=\mu(\rho^{-1}A)$ for all $A$ in the domain of $\mu$.

\begin{lemma}\label{L1}Let $L,E\in G(d,k)$, $\rho\in \mathit{SO}_d$ and
$L=\rho E$. If $|\rho|\le1/8$, then
\[
d_P(\pi_L\varphi,\rho\pi_E\varphi)\le3|\rho|^{1/3}.
\]
\end{lemma}

\begin{pf}Put $\mu:=\pi_L\varphi$ and $\nu:= \rho\pi_E\varphi
$, write $\varepsilon:=|\rho|^{1/3}$.

We have to show that $\mu(A)\le\nu(A_{3\varepsilon})+3\varepsilon$ and
$\nu(A)\le\mu(A_{3\varepsilon})+3\varepsilon$ for all $A\in\mathcal{
B}({\mathbb S}_L^{k-1})$. Let $A\in\mathcal{B}({\mathbb S}_L^{k-1})$ be
given. The first assertion reads
\begin{eqnarray}\label{4}
&& \int_{{\mathbb S}^{d-1}\setminus L^\perp} \mathbf{1}_{A}({\rm
pr}_L(\mathbf{u}))\|\mathbf{u}|L\| \varphi(\mathrm{d}\mathbf{u})\nonumber
\\[-8pt]\\[-8pt]
&&\qquad \le\int_{{\mathbb S}^{d-1}\setminus E^\perp} \mathbf{1}_{\rho
^{-1}A_{3\varepsilon}}(\operatorname{pr}_E(\mathbf{u}))\|\mathbf{u}|E\| \varphi
(\mathrm{d}\mathbf{u}
)+3\varepsilon.\nonumber
\end{eqnarray}
Writing
\begin{eqnarray*}
M_1 &=& \{\mathbf{u}\in{\mathbb S}^{d-1}\dvtx \|\mathbf{u}|L\|<\varepsilon
\},\\
M_2 &=& \{\mathbf{u}\in{\mathbb S}^{d-1}\dvtx \|\mathbf{u}|L\|\ge
\varepsilon\},
\end{eqnarray*}
we have
\[
\int_{M_1\setminus L^\perp} \mathbf{1}_{A}(\operatorname{pr}_L(\mathbf{u}))\|
\mathbf{u}|L\|
\varphi(\mathrm{d}\mathbf{u}) \le\varepsilon.
\]

Let $\mathbf{u},\mathbf{v}\in{\mathbb S}^{d-1}$ and assume that $\|
\mathbf{u}-\mathbf{v}\|\le
\varepsilon^3$ and $\mathbf{u}\in M_2$, hence $ \|\mathbf{u}|L\|\ge
\varepsilon$. From
\[
\bigl| \|\mathbf{u}|L\|-\|\mathbf{v}|L\|\bigr| \le\|(\mathbf
{u}-\mathbf{v})|L\| \le\|\mathbf{u}-\mathbf{v}
\| \le\varepsilon^3
\]
we get $ \|\mathbf{v}|L\| \ge\|\mathbf{u}|L\| -\varepsilon^3 \ge
\varepsilon-\varepsilon
^3 \ge\varepsilon/2>0$, hence $\mathbf{v}\in{\mathbb S}^{d-1}\setminus
L^\perp$.

There are unique representations
\begin{eqnarray*}
\mathbf{u}&=& t\mathbf{u}_0+\mathbf{u}_1, \qquad\mathbf{u}_0\in
L\cap{\mathbb S}^{d-1},  \mathbf{u}
_1\in L^\perp, t>0,
\\
\mathbf{v}&=& \tau\mathbf{v}_0+\mathbf{v}_1, \qquad\mathbf
{v}_0\in L\cap{\mathbb S}^{d-1},
\mathbf{v}_1\in L^\perp, \tau>0.
\end{eqnarray*}
Here, $\mathbf{u}_0=\operatorname{pr}_L(\mathbf{u})$ and $\mathbf{v}_0={\rm
pr}_L(\mathbf{v})$.
From $|t-\tau| = |\|\mathbf{u}|L\|-\|\mathbf{v}|L\|| \le
\varepsilon^3$
together with $\tau=\|\mathbf{v}|L\|\ge\varepsilon/2$ and $t=\|\mathbf
{u}|L\|\ge
\varepsilon$, we get
\begin{eqnarray*}
\|\mathbf{u}_0-\mathbf{v}_0\| &\le& \biggl\|\frac{\mathbf
{u}}{t}-\frac{\mathbf{v}}{\tau
}\biggr\| =\biggl\|\frac{\mathbf{u}}{t}-\frac{\mathbf{v}}{t}+\frac
{\mathbf{v}
}{t}-\frac{\mathbf{v}}{\tau}\biggr\|
\\
&\le& \frac{1}{t}\|\mathbf{u}-\mathbf{v}\| + \biggl|\frac
{1}{t}-\frac{1}{\tau
}\biggr| \le\varepsilon^{-1}\cdot\varepsilon^3(1+2/\varepsilon) \le
3\varepsilon.
\end{eqnarray*}
Now let $\mathbf{v}:=\rho\mathbf{u}$, then $\|\mathbf{u}-\mathbf
{v}\|\le|\rho|=\varepsilon^3$.
Hence, we have $\|\mathbf{u}_0-\mathbf{v}_0\|\le3\varepsilon$.

Suppose that $\mathbf{u}\in M_2$ is such that $\mathbf{u}_0\in A$.
Then ${\rm
pr}_L(\rho\mathbf{u}) = \mathbf{v}_0 \in A_{3\varepsilon}$, hence
$\operatorname{pr}_E(\mathbf{u})
\in\rho^{-1}A_{3\varepsilon}$. From
\[
\|\mathbf{u}|L\| - \|\mathbf{u}|E\| = \|\mathbf{u}|L\|-\|\rho
\mathbf{u}|L\| = \|\mathbf{u}|L\|-\|\mathbf{v}
|L\|\le\varepsilon^3
\]

\noindent
we see that $\mathbf{u}\notin E^\perp$ and conclude that
\begin{eqnarray*}
&& \int_{M_2} \mathbf{1}_A(\operatorname{pr}_L(\mathbf{u}))\|\mathbf{u}|L\|
\varphi(\mathrm{d}\mathbf{u})
\\
&&\qquad \le\int_{{\mathbb S}^{d-1}\setminus E^\perp} \mathbf{1}_{\rho
^{-1}A_{3\varepsilon}}(\operatorname{pr}_E(\mathbf{u}))\|\mathbf{u}|L\| \varphi
(\mathrm{d}\mathbf{u})\\
&&\qquad \le\int_{{\mathbb S}^{d-1}\setminus E^\perp} \mathbf{1}_{\rho
^{-1}A_{3\varepsilon}}(\operatorname{pr}_E(\mathbf{u}))\|\mathbf{u}|E\| \varphi
(\mathrm{d}\mathbf{u})
+\varepsilon^3.
\end{eqnarray*}
Altogether, we obtain
\begin{eqnarray*}
&& \int_{{\mathbb S}^{d-1}\setminus L^\perp} \mathbf{1}_{A}({\rm
pr}_L(\mathbf{u}))\|\mathbf{u}|L\| \varphi(\mathrm{d}\mathbf{u})
\\
&& \qquad\le\varepsilon+\int_{{\mathbb S}^{d-1}\setminus E^\perp} \mathbf{
1}_{\rho^{-1}A_{3\varepsilon}}(\operatorname{pr}_E(\mathbf{u}))\|\mathbf{u}|E\|
 \varphi(\mathrm{d}
\mathbf{u})+\varepsilon^3
\end{eqnarray*}
and hence (\ref{4}).

Since $|\rho|=|\rho^{-1}|$, inequality (\ref{4}) remains true if we
interchange $L$ with $E$ and $\rho$ with $\rho^{-1}$ and then replace
$A$ by $\rho^{-1}A$. The resulting inequality is $\nu(A) \le\mu
(A_{3\varepsilon}) + 3\varepsilon$, which completes the proof of Lemma \ref{L1}.
\end{pf}

In the following, the dependence of some constants $c_i$ on the measure
$\varphi$ is only via the number
\[
m(\varphi):= \min_{\mathbf{u}\in{\mathbb S}^{d-1}} \int_{{\mathbb
S}^{d-1}} |\langle\mathbf{u},\mathbf{v}\rangle| \varphi(\mathrm
{d}\mathbf{v}).
\]
This number, which can be considered as a measure of nondegeneracy, is
positive, since the support of $\varphi$ is not contained in a great subsphere.

\begin{lemma}\label{L2}Let $B$ be the $\mathbf{o}$-symmetric convex
body with
$S_{d-1}(B,\cdot) =\varphi$. The inradius $r$ and circumradius $R$ of
$B$ can be estimated by
\[
c_4 \le r\le R\le c_5,
\]
where $c_4,c_5>0$ are constants depending only on $m(\varphi)$ and an
upper bound for~$\varphi$. (Here, $\varphi$ can be any finite even
measure on ${\mathbb S}^{d-1}$ not concentrated on a great subsphere.)
\end{lemma}

\begin{pf} First, we repeat a known argument (\cite{Sch93}, page 303).
If $\varphi({\mathbb S}^{d-1})\le b$, the isoperimetric inequality
gives $V_d(B)\le c(b)$, with a constant $c(b)$ depending only on $b$
(and the dimension, which we do not mention). Let $\mathbf{x}\in B$. Then
\begin{eqnarray*}
V_d(B) &=& \frac{1}{d} \int_{{\mathbb S}^{d-1}} h(B,\mathbf{v})
\varphi
(\mathrm{d}\mathbf{v}) \ge\frac{1}{d} \int_{{\mathbb S}^{d-1}}
\max\{\langle\mathbf{x}
,\mathbf{v}\rangle,0\}   \varphi(\mathrm{d}\mathbf{v})
\\
&=& \frac{1}{2d} \int_{{\mathbb S}^{d-1}} |\langle\mathbf
{x},\mathbf{v}\rangle|
  \varphi(\mathrm{d}\mathbf{v}) \ge\frac{1}{2d}\|\mathbf{x}\|
m(\varphi).
\end{eqnarray*}
It follows that $R\le2dc(b)/m(\varphi)$.

Second, since $B$ is centrally symmetric, an inball of $B$ is touched
by two parallel supporting hyperplanes of $B$. Let $\mathbf{u}$ be a unit
vector parallel to these hyperplanes. The projection $B|\mathbf
{u}^\perp$
lies between two parallel hyperplanes in $\mathbf{u}^\perp$ which are
distance $2r$ apart, and it is contained in $R{\mathbb B}^d\cap\mathbf{u}
^\perp$. Hence, $V_{d-1}(B|\mathbf{u}^\perp)\le2rV_{d-2}({\mathbb
B}^{d-2})R^{d-2}$. On the other hand, using \cite{Sch93}, (7.4.1),
\begin{eqnarray*}
V_{d-1}(B|\mathbf{u}^\perp)
&=& \frac{1}{2} \int_{{\mathbb S}^{d-1}} |\langle\mathbf{u},\mathbf
{v}\rangle|
S_{d-1}(B,\mathrm{d}\mathbf{v})
\\
&=& \frac{1}{2} \int_{{\mathbb S}^{d-1}} |\langle\mathbf{u},\mathbf
{v}\rangle|
\varphi(\mathrm{d}\mathbf{v}) \ge\frac{1}{2}m(\varphi).
\end{eqnarray*}
The assertion follows.
\end{pf}

It is technically convenient to consider also the $\mathbf
{o}$-symmetric convex
body $B(L)$ in $L$ with surface area measure
\[
S^L_{k-1}(B(L),\cdot) = \pi_L\varphi.
\]
From (\ref{5a}), we have
%
\begin{equation}\label{4.3a}
B_L=\biggl(\frac{\gamma}{\gamma_{X\cap L}}\biggr)^{{1/(k-1)}}B(L)
\end{equation}
and
%
\begin{equation}\label{4.3}
\frac{\gamma_{ X\cap L}}{\gamma}= \int_{{\mathbb S}^{d-1}} \|
\mathbf{v}
|L\| \varphi(\mathrm{d}\mathbf{v}).
\end{equation}

\begin{lemma}\label{L3}Let $k\in\{2,\dots,d-1\}$. Let $r_L,R_L$
denote the inradius and circumradius, respectively, of either $B(L)$ or
$B_L$, measured in $L\in G(d,k)$. There are constants $c_6,c_7>0$,
depending only on $m(\varphi)$, such that
\[
c_6\le r_L\le R_L\le c_7\qquad\mbox{for all }L\in G(d,k).
\]
\end{lemma}

\begin{pf}Let $L\in G(d,k)$. From (\ref{0}), clearly $\pi_L\varphi
({\mathbb S}_L^{k-1}) \le\varphi({\mathbb S}^{d-1})$, hence $1$ is an
upper bound for $\pi_L\varphi$.

Let $\Pi_{\varphi/2}$ be the zonoid with generating measure $\varphi
/2$. Then $\Pi_{\varphi/2}|L$ has generating measure $\pi_L\varphi
/2$. For $\mathbf{u}\in L$, it follows that
\begin{eqnarray*}
\int_{{\mathbb S}_L^{k-1}} |\langle\mathbf{u},\mathbf{v}\rangle|
\pi_L\varphi
(\mathrm{d}\mathbf{v}) &=& 2 h(\Pi_{\varphi/2}|L,\mathbf{u})
\\
& = & 2h(\Pi_{\varphi/2},\mathbf{u}) = \int_{{\mathbb S}^{d-1}}
|\langle
\mathbf{u},\mathbf{v}\rangle| \varphi(\mathrm{d}\mathbf{v})\ge
m(\varphi).
\end{eqnarray*}
Now Lemma \ref{L2}, applied in $L$ and to the measure $\pi_L\varphi
$, shows that the inradius and circumradius of $B(L)$ can be estimated
from both sides by positive constants depending only on $m(\varphi)$.
The same fact for $B_L$ follows from (\ref{4.3a}), since (\ref{4.3}) gives
%
\begin{equation}\label{spaet2}
m(\varphi) \le\gamma_{X\cap L}/\gamma\le1.
\end{equation}
For the left side, note that $L$ contains a unit vector $\mathbf{u}$
and that
$\|\mathbf{v}|L\|\ge|\langle\mathbf{v},\mathbf{u}\rangle|$.
\end{pf}

In the following, we make use of the Hausdorff metric $\delta$ and of
the deviation function $\vartheta$ defined by (\ref{n2}).

\begin{lemma}\label{L4} Let $k\in\{2,\dots,d-1\}$. There exist
constants $c_{8}, c_{9}$, depending only on $m(\varphi)$, with the
following property. If $L,E\in G(d,k)$ and if $\rho\in \mathit{SO}_d$ is a
rotation with $L=\rho E$ and $|\rho|\le1/8$, then
%
\begin{equation}\label{4a}
\delta(B(L),\rho B(E)) \le c_{8}|\rho|^{1/3k}
\end{equation}
and
%
\begin{equation}\label{4d}
\vartheta(B_L,\rho B_E) \le c_{9}|\rho|^{1/3k}.
\end{equation}
\end{lemma}

\begin{pf}By Lemma \ref{L3}, the inradius and circumradius of
$B(L')$, $L'\in G(d,k)$, can be bounded from below and from above by
positive constants depending only on $m(\varphi)$. We use the
stability result of \cite{HS02}, Theorem 3.1, for the solutions of
Minkowski's problem and apply it here in the subspace $L$ of the
assertion. [Note that $B(L),B(E)$ are centrally symmetric, hence the
translations appearing {\em loc.~cit.} can be omitted.] We conclude that
\begin{eqnarray*}
\delta(B(L),\rho B(E)) &\le& c  d_P(S^L_{k-1}(B(L),\cdot),
S^L_{k-1}(\rho B(E),\cdot))^{1/k}
\\
&=& c  d_P(\pi_L\varphi,\rho\pi_E\varphi)^{1/k}
\end{eqnarray*}
with a constant $c$ depending only on $m(\varphi)$. By Lemma \ref{L1},
\[
d_P(\pi_L\varphi,\rho\pi_E\varphi) \le3|\rho|^{1/3},
\]
hence (\ref{4a}) follows, with suitable $c_{8}$.

From (\ref{4a}), with $\lambda:= c_{8}|\rho|^{1/3k}$, we get
$B(L)\subset\rho B(E)+\lambda{\mathbb B}^k_L$, where ${\mathbb
B}^k_L:={\mathbb B}^d\cap L$ is the unit ball in $L$. Since $c_6
{\mathbb B}^k_L\subset\rho B(E)$ by Lemma \ref{L3}, we get
\[
B(L) \subset(1+\lambda/c_6)\rho B(E).
\]
A similar relation holds with $\rho B(E)$ and $B(L)$ interchanged,
hence
\[
(1+\lambda/c_6)^{-1}\rho B(E) \subset B(L) \subset(1+\lambda
/c_6)\rho B(E).
\]
This gives
\[
\vartheta(B(L),\rho B(E)) \le\log(1+\lambda/c_6)^2 \le2\lambda
/c_6.
\]
Here, we may replace $B(L),B(E)$ by $B_L,B_E$, since $\vartheta$ is
invariant under dilatations.
\end{pf}

In the following lemma, the deviation function $\vartheta$ for convex
bodies in different subspaces, as defined by (\ref{n3}), is used.

\begin{lemma}\label{L5}Let $k\in\{2,\dots,d-1\}$, let $L,L^*\in
G(d,k)$ and $\varepsilon>0$. If $\Delta(L,L^*) \le\min\{1/8,(\varepsilon
/c_{9})^{3k}\}$, where $c_{9}$ is the constant appearing in Lemma \ref
{L4}, then every convex body $K\in\mathcal{K}_0(L)$ with $\vartheta
(K,B_L) < \varepsilon$ satisfies $\vartheta(K,B_{L^*}) < 2\varepsilon$.
\end{lemma}

\begin{pf}Let $L,L^*\in G(d,k)$ and $\Delta(L,L^*) \le\min\{
1/8,(\varepsilon/c_{9})^{3k}\}$. There exists a rotation $\rho\in \mathit{SO}_d$
with $\rho L=L^*$ and $|\rho|=\Delta(L,L^*)$, hence $\rho$ satisfies
$|\rho|\le(\varepsilon/c_{9})^{3k}$ and $|\rho|\le1/8$. Lemma \ref
{L4} gives
%
\begin{equation}\label{2a}
\vartheta(B_{L^*},\rho B_L)\le\varepsilon.
\end{equation}

Suppose that $K\in\mathcal{K}_0(L)$ and $\vartheta(K,B_L) < \varepsilon$.
The definition of $\vartheta$ implies the triangle inequality
\[
\vartheta(\rho K,B_{L^*}) \le\vartheta(\rho K,\rho B_L) +\vartheta
(\rho B_L,B_{L^*}).
\]
In fact, if $\vartheta(\rho K,\rho B_L)=:a_1$ and $\vartheta(\rho
B_L,B_{L^*})=:a_2$, there are numbers $\alpha_1,\beta_1,\break \alpha
_2,\beta_2>0$ with
$\log(\beta_1/\alpha_1) =a_1$, $\log(\beta_2/\alpha_2)=a_2$ and a
vector $\mathbf{z}\in\rho L$ such that
\[
\alpha_1 \rho B_L \subset\rho K+\mathbf{z}\subset\beta_1 \rho
B_L,\qquad
\alpha_2 B_{L^*} \subset\rho B_L \subset\beta_2 B_{L^*}.
\]
In the second case, we have used that, due to the central symmetry of
$B_L$ and $B_{L^*}$, the translation vector appearing in the definition
of $\vartheta$ can be omitted. We deduce that
\[
\alpha_1\alpha_2 B_{L^*} \subset\rho K+\mathbf{z}\subset\beta
_1\beta_2 B_{L^*}
\]
and hence $\vartheta(\rho K,B_{L^*}) \le\log(\beta_1\beta_2/\alpha
_1\alpha_2)= a_1+a_2$.

From $\vartheta(\rho K,\rho B_L)=\vartheta(K,B_L)<\varepsilon$ and
$\vartheta(\rho B_L,B_{L^*}) = \vartheta(B_{L^*},\rho B_L) \le
\varepsilon$ we get $\vartheta(\rho K, B_{L^*}) <2\varepsilon$. Since
$|\rho|= \Delta(L,L^*)$, this yields $\vartheta(K,B_{L^*}) <
2\varepsilon$, which finishes the proof.
\end{pf}

\section{Two preparatory probability estimates}\label{sec4}

The plan is to prove Theorem~\ref{T1} by using (\ref{n1}) and
applying results for the zero cell $Z_0$ of $X$, obtained in \cite
{HRS04a}, to the random polytope $Z_0\cap L$, for each $L\in G(d,k)$.
This is possible since $Z_0\cap L$ is stochastically equivalent to the
zero cell of the section
process $X\cap L$, which is a stationary Poisson hyperplane process in
$L$; it has intensity $\gamma_{X\cap L}$ and spherical directional
distribution $\varphi_{X\cap L} = (\gamma/\gamma_{X\cap L}) \pi
_L\varphi$, by (\ref{5a}). In applying the results from \cite
{HRS04a}, we have to ensure that the constants appearing there can be
chosen independently of $L$.

Let $L\in G(d,k)$. For $\mathbf{u}\in{\mathbb S}^{k-1}_L$ and $t>0$
we write
$H_L^-(\mathbf{u},t):=H^-(\mathbf{u},t)\cap L$. If $\mathbf
{u}_1,\dots,\mathbf{u}_n\in
{\mathbb S}^{k-1}_L$ and $t_1,\dots,t_n>0$, we use the notation
\[
\bigcap_{i=1}^n H_L^-(\mathbf{u}_i,t_i)=:P_L(\mathbf{u}_1,\ldots
,\mathbf{u}_n;t_1,\ldots,t_n).
\]
In the following, we write
\[
\tau_L:= kV_k(B_L)^{1-1/k}.
\]

\begin{lemma}\label{L6} Let $\beta>0$. There are positive constants
$c_{10},h_0$, depending only on $\varphi, \gamma$ and $\beta$, such
that for all $L \in G(d,k)$, all $a\ge1$ and $0<h\le h_0$,
\[
\mathbf{P}\{V_k(Z_0\cap L) \in a(1,1+h)\} \ge c_{10} h \exp
[-2(1+\beta
)\gamma_{X\cap L}\tau_L a^{1/k}].
\]
\end{lemma}

\begin{pf} Let $\beta>0$ and $a\ge1$ be given.

First, we consider a fixed $L^*\in G(d,k)$. Then $X\cap L^*$ is a stationary
Poisson hyperplane process
in $L^*$ with intensity $\gamma_{X\cap L^*}$ and spherical directional
distribution $\varphi_{X\cap L^*}$. For given $\beta>0$, Lemma 3.1 in
\cite{HRS04a}, applied to the convex body $B_{L^*}$ in $L^*$, yields
the existence of a number $N\in\mathbb{N}$, of unit vectors $\mathbf
{u}_1^0,\dots
,\mathbf{u}_N^0
\in{\mathbb S}_{L^*}^{k-1}$ in the support of the measure $\varphi
_{X\cap L^*}$ (which is equal to the support of $\pi_{L^*}\varphi$)
and of positive numbers $t_1^0,\dots,t_N^0$, all depending only on
$\varphi, \gamma, L^*$ and $\beta$, such that the polytope
\[
P^0:= P_{L^*}(\mathbf{u}_1^0,\ldots,\mathbf{u}_N^0;t_1^0,\ldots,t_N^0)
\]
has $N$ facets (in $L^*$) and satisfies
\[
P^0\subset(1+\beta/4)B_{L^*} \quad\mbox{and}\quad V_k(P^0)=V_k(B_{L^*}).
\]
Next, we can choose neighbourhoods $U_i$ of $\mathbf{u}_i^0$ in
${\mathbb
S}^{k-1}_{L^*}$ and a number $\alpha>0$ with $t_i^0-\alpha>0$,
$i=1,\ldots,N$, all depending only on $\varphi, \gamma, L^*$ and
$\beta$, such that, for all $\mathbf{u}_1,\ldots,\mathbf{u}_N\in
{\mathbb
S}^{k-1}_{L^*}$ and $t_1,\ldots,t_N\in\mathbb{R}$
with
%
\begin{equation}\label{4.1}
\mathbf{u}_i\in U_i,\qquad|t_i-t_i^0|<\alpha,\qquad i=1,\ldots,N,
\end{equation}
the following condition (i) is satisfied.
\begin{longlist}
\item[(i)] $P:= P_{L^*}(\mathbf{u}_1,\dots,\mathbf{u}_N;t_1,\dots,t_N)$ is
a polytope in
$L^*$ with $N$ facets and satisfying $P\subset(1+\beta/2)B_{L^*}$.
\end{longlist}

The set of values
\[
V_k(P_{L^*}(\mathbf{u}_1^0,\ldots,\mathbf{u}_N^0;t_1^0,\ldots
,t_{N-1}^0,t)) \qquad
\mbox{with }|t-t_N^0|<\alpha
\]
is an interval containing $V_k(B_{L^*})$ in its interior. Therefore,
after decreasing $U_1,\ldots,U_N,\alpha$, if necessary, we can assume
that there exists a number $b>0$, depending only on $\varphi, \gamma,
L^*$ and $\beta$, with the following property.
\begin{longlist}
\item[(ii)] If (\ref{4.1}) is satisfied, then
\begin{eqnarray*}
&& \bigl(V_k(B_{L^*})-b,V_k(B_{L^*})+b\bigr)
 \\
&& \qquad\subset\{V_k(P_{L^*}(\mathbf{u}_1,\ldots,\mathbf{u}_N;t_1,\dots
,t_{N-1},t))\dvtx
|t-t_N^0|<\alpha\}.
\end{eqnarray*}
\end{longlist}

We must extend the preceding to the subspaces $L$ in a suitable
neighbourhood $N_\theta(L^*)$. The numbers $\theta, \eta, h_0$
appearing in the following can be chosen to depend only on $\varphi,
\gamma, L^*$ and $\beta$. Let $\theta\in(0,1/8]$; below it will be
specified further. To each $L\in N_\theta(L^*)$, we choose a rotation
$\rho_L$ with $L=\rho_L L^*$ and $|\rho_L|\le\theta$.

We choose a number $\eta>0$ so small that to each $i\in\{1,\dots,N\}
$ there exists a neighbourhood $U_i'$ of $\mathbf{u}_i^0$ in ${\mathbb
S}^{k-1}_{L^*}$ with $(U'_i)_\eta\subset U_i$, where for $A\subset
{\mathbb S}^{k-1}_{L^*}$ the set $A_\eta$ is defined by $A_\eta=\{
\mathbf{y}
\in{\mathbb S}^{k-1}_{L^*}\dvtx \|\mathbf{y}-\mathbf{x}\|<\eta\mbox{
for some }\mathbf{x}\in
A\}$. Decreasing $\eta$, if necessary (without changing the sets
$U'_i$), we can also assume that
\[
\pi_{L^*}\varphi(U_i') \ge2\eta\qquad\mbox{for }i=1,\ldots,N.
\]
This is possible since $U_i'$ is a neighbourhood of $\mathbf{u}_i^0\in
\operatorname{supp} \pi_{L^*}\varphi$.

By Lemma \ref{L1}, we can further choose $\theta$ so small that
\[
d_P(\pi_L\varphi,\rho_L\pi_{L^*}\varphi)\le\eta\qquad\mbox{for
}L\in N_\theta(L^*).
\]
Then,
\[
\pi_L\varphi(\rho_L U_i) \ge\pi_L\varphi((\rho_L U_i')_\eta)
\ge(\rho_L \pi_{L^*}\varphi)(\rho_L U_i') -\eta= \pi
_{L^*}\varphi(U_i') -\eta\ge\eta.
\]
Hence, putting $ U_i^L:= \rho_L U_i$, we have from (\ref{5a})
\[
\varphi_{X\cap L}(U_i^L) =\frac{\gamma}{\gamma_{X\cap L}}\pi
_L\varphi(U_i^L) \ge\eta>0 \qquad\mbox{for } i=1,\dots,N.
\]
Due to (\ref{4d}), we can decrease $\theta$, if necessary, such that
%
\begin{equation}\label{eqn1}
\rho_L B_{L^*} \subset\frac{1+\beta}{1+\beta/2}B_L \qquad\mbox
{for } L\in N_\theta(L^*).
\end{equation}
Using (ii) above, \eqref{4a} and the fact that $L\mapsto\gamma
_{X\cap L}=\gamma\pi_L\varphi(\mathbb{S}^{k-1}_L)$ is continuous by
Lemma \ref{L1}, we can decrease $\theta$ further, if necessary, and
choose a number $h_0>0$ such that
%
\begin{equation}\label{4.5}
L\in N_\theta(L^*), \qquad\mathbf{u}_i\in U^L_i, \qquad|t_i-t_i^0|
<\alpha
,\qquad i=1,\dots,N,
\end{equation}
implies
\[
V_k(B_L)(1,1+h_0) \subset\{V_k(P_L(\mathbf{u}_1,\dots,\mathbf
{u}_N;t_1,\ldots
,t_{N-1},t))\dvtx |t-t_N^0|<\alpha\}.
\]
Here, we have used that
\[
P_L(\mathbf{u}_1,\ldots,\mathbf{u}_N;t_1,\ldots,t_{N-1},t) = \rho_L
P_{L^*}(\rho
_L^{-1}\mathbf{u}_1,\ldots,\rho_L^{-1}\mathbf{u}_N;t_1,\ldots,t_{N-1},t).
\]

After these choices, the following is true for all $L\in N_\theta
(L^*)$. If (\ref{4.5}) holds, then (${\rm i}_L$) and (${\rm ii}_L$)
are satisfied:
\begin{longlist}[(${\rm ii}_L$)]
\item[(${\rm i}_L$)] $P:= P_L(\mathbf{u}_1,\ldots,\mathbf{u}_N;t_1,\ldots
,t_N)$ is a
polytope with $N$ facets and satisfying
$P\subset(1+\beta)B_L$.\vspace*{8pt}

\item[(${\rm ii}_L$)] \[
 V_k(B_L)(1,1+h_0) \subset\{
V_k(P_L(\mathbf{u}_1,\ldots,\mathbf{u}_N;t_1,\ldots,t_{N-1},t))\dvtx
|t-t_N^0|<\alpha\}.\]
\end{longlist}
 In fact, (${\rm i}_L$) follows from (${\rm i}$) and \eqref
{eqn1}, since $\rho_L^{-1}\mathbf{u}_i\in U_i$ and therefore
\begin{eqnarray*}
P&=&\rho_LP_{L^*}(\rho_L^{-1}\mathbf{u}_1,\ldots,\rho
_L^{-1}\mathbf{u}
_N;t_1,\ldots,t_{N-1},t)
\\
&\subset&\rho_L(1+\beta/2)B_{L^*}\subset(1+\beta
/2)\frac{1+\beta}{1+\beta/2}B_L
=(1+\beta)B_L.
\end{eqnarray*}

We restate what we have found so far, making explicit the dependence on
$L^*$. For any $L^*\in G(d,k)$, there exist numbers $\theta(L^*)\in
(0,1/8]$, $N(L^*)\in\mathbb{N}$, $\alpha(L^*)>0$, $t_1^0(L^*),\ldots
,t^0_{N(L^*)}(L^*)>\alpha(L^*)$, $h_0(L^*)>0$, $\eta(L^*)>0$, unit vectors
$\mathbf{u}^0_1(L^*),\ldots,\mathbf{u}^0_{N(L^*)}(L^*)\in\mathbb
{S}_{L^*}^{k-1}$
and neighbourhoods $U_i(L^*)$ of $\mathbf{u}^0_i(L^*)$ in $\mathbb
{S}^{k-1}_{L^*}$, $i=1,\ldots,N(L^*)$, such that for all $L\in
N_{\theta(L^*)}(L^*)$ and for $\mathbf{u}_i\in U_i^L(L^*)$ and
$|t_i-t_i^0(L^*)|<\alpha(L^*)$, $i=1,\ldots,N(L^*)$, the following
conditions are satisfied:
\[
\varphi_{X\cap L}(U_i^L(L^*))\ge\eta(L^*)>0,\qquad i=1,\ldots,N(L^*),
\]\vspace*{-15pt}
\begin{longlist}[(${\rm ii}_L$)]
\item[(${\rm i}_L$)] $P:= P_L(\mathbf{u}_1,\ldots,\mathbf
{u}_{N(L^*)};t_1,\ldots
,t_{N(L^*)})$ is a polytope with $N(L^*)$ facets and satisfying
$P\subset(1+\beta)B_L$, and

\item[(${\rm ii}_L$)]
\begin{eqnarray*}
&&V_k(B_L)\bigl(1,1+h_0(L^*)\bigr)
\\
&&\qquad \subset\bigl\{
V_k\bigl(P_L\bigl(\mathbf{u}_1,\ldots,\mathbf{u}_{N(L^*)};t_1,\ldots,t_{N(L^*)-1},t\bigr)\bigr)\dvtx
\bigl|t-t_{N(L^*)}^0\bigr|<\alpha(L^*)\bigr\}.
\end{eqnarray*}
\end{longlist}

Since $(G(d,k),\Delta)$ is compact and $\{N_{\theta(L^*)}(L^*):
L^*\in G(d,k)\}$ is an open cover of $G(d,k)$, there are $L_1^*,\ldots
,L_r^*\in G(d,k)$ such that $\{N_{\theta(L^*_j)}(L^*_j)\dvtx j=1,\ldots,r\}
$ is a finite subcover of $G(d,k)$. We put
\begin{eqnarray*}
\eta_0&:=&\min\{\eta(L_j^*):j=1,\ldots,r\}>0,
\\
h_0&:=&\min\{h_0(L_j^*):j=1,\ldots,r\}>0.
\end{eqnarray*}
Hence, for $L\in G(d,k)$ there is some $j\in\{1,\ldots,r\}$ such that
$L\in N_{\theta(L^*_j)}(L^*_j)$ and
%
\begin{equation}\label{eqnew}
\varphi_{X\cap L}(U_i^L(L^*_j))\ge\eta(L^*_j)\ge\eta_0>0.
\end{equation}
Note that $U_i^L(L_j^*)=\rho_L(L_j^*)U_i(L_j^*)$. For
$\mathbf{u}
_i\in U_i^L(L_j^*)$ and for $|t_i-t_i^0(L^*_j)|<\alpha(L^*_j)$,
$i=1,\ldots,N(L_j^*)$, the set
\[
P:= P_L\bigl(\mathbf{u}_1,\ldots,\mathbf{u}_{N(L^*_j)};t_1,\ldots,t_{N(L^*_j)}\bigr)
\]
is a polytope with $N(L^*_j)$ facets and satisfying $P\subset(1+\beta
)B_L$, and
\begin{eqnarray*}
&&V_k(B_L)(1,1+h_0)
\\
&&\qquad \subset\bigl\{V_k\bigl(P_L\bigl(\mathbf{u}_1,\ldots,\mathbf
{u}_{N(L^*_j)};t_1,\ldots
,t_{N(L^*_j)-1},t\bigr)\bigr)\dvtx \bigl|t-t_{N(L^*_j)}^0\bigr|<\alpha(L^*_j)\bigr\}.
\end{eqnarray*}

We are now in a situation where we can adjust the second part of the
proof of
\cite{HRS04a}, Lemma 3.2, in a fixed linear subspace $L\in G(d,k)$. We
choose a
corresponding index $j\in\{1,\ldots,r\}$ such that $L\in N_{\theta
(L^*_j)}(L^*_j)$,
as described above. For the given $a\ge1$, we choose a number $\varrho
>0$ such that $V_k(\varrho B_L)=a$,
that is, $\varrho=a^{1/k}V_k(B_L)^{-1/k}$. Then, for $\mathbf{u}_i\in
U_i^L(L^*_j)$ and for
$|t_i-t_i^0(L_j^*)|<\alpha(L_j^*)$, $i=1,\ldots,N(L_j^*)$,
\begin{longlist}[(${\rm ii}_\varrho$)]
\item[(${\rm i}_\varrho$)] $P_\varrho:= P_L(\mathbf{u}_1,\ldots,\mathbf{u}
_{N(L_j^*)};\varrho t_1,\ldots,\varrho
t_{N(L_j^*)})$
is a polytope with $N(L^*_j)$ facets and satisfying $P_\varrho\subset
(1+\beta)\varrho B_L$, and,

\item[(${\rm ii}_\varrho$)] for $0<h\le h_0$,
\[
V_k(B_L)(1,1+h)\subset\bigl\{v_L(t)\dvtx\bigl|t-\varrho
t_{N(L_j^*)}^0(L_j^*)\bigr|<\varrho\alpha(L_j^*)\bigr\}
\]
with
\[
v_L(t):=V_k\bigl(P_L\bigl(\mathbf{u}_1,\ldots,\mathbf{u}_{N(L_j^*)};\varrho
t_1,\ldots
,\varrho t_{N(L_j^*)-1},t\bigr)\bigr).
\]
\end{longlist}
Let $\lambda^1$ denote 1-dimensional Lebesgue measure.
The argument on page 1147, lines~$-$17 to bottom, in \cite{HRS04a} now
shows that
\begin{eqnarray*}
&&\lambda^1\bigl\{t\in\mathbb{R}\dvtx \bigl|t-\varrho
t_{N(L_j^*)}^0(L_j^*)\bigr|<\varrho
\alpha(L_j^*),
v_L(t)\in V_k(B_L)(1,1+h)\bigr\}
\\
&&\qquad\ge c(\beta,\varphi)\varrho h,
\end{eqnarray*}
where $c(\beta,\varphi)$ is a constant depending only on $\beta$ and
$\varphi$.
Here it is implicitly used that $P_\varrho\subset(1+\beta)\varrho
B_L$, which implies that the $(k-1)$-dimensional
volume of the orthogonal projection of $P_\varrho$ on to the
orthogonal complement of $\mathbf{u}_{N(L_j^*)}$
can be bounded from above by a constant depending only on $\beta$ and
$m(\varphi)$. Moreover, it is also
used that $V_k(B_L)$ can be bounded from below by a constant depending
only on $m(\varphi)$.

Next, we define a sufficiently large set of convex polytopes in $L$ by
\begin{eqnarray*}
\mathcal{P}_L &:=& \bigl\{P_L\bigl(\mathbf{u}_1,\ldots,\mathbf{u}
_{N(L_j^*)};t_1,\ldots,t_{N(L_j^*)}\bigr)\dvtx \mathbf{u}_i\in U_i^L(L_j^*)
\mbox{ and }
\\
& &\ \ |t_i-\varrho t_i^0(L_j^*)|<\varrho\alpha(L_j^*), \mbox{ for
}i=1,\ldots,N(L_j^*), \mbox{ and }
\\
& &\ \ V_k\bigl(P_L\bigl(\mathbf{u}_1,\ldots,\mathbf{u}_{N(L_j^*)};t_1,\ldots
,t_{N(L_j^*)}\bigr)\bigr)\subset V_k(\varrho B_L)(1,1+h)\bigr\}.
\end{eqnarray*}

Let $\mathcal{H}_{(1+\beta)\varrho B_L}:=\{H\in A(d,k-1)\dvtx(1+\beta
)\varrho B_L\cap H\neq\varnothing\}$.
For a hyperplane process $Y$ in $L$, we write $Z_0(Y)$ for the induced
zero cell in $L$, and ``$\,\rule{.1mm}{.26cm}\rule{.24cm}{.1mm}\, $'' denotes the restriction of a measure.
Subsequently, we adapt the argument from \cite{HRS04a}, page~1148,
to the present situation. For the first estimate, we use that any
polytope in $\mathcal{P}_L$
is contained in $(1+\beta)\varrho B_L$. Thus, we get
\begin{eqnarray*}
&&\mathbf{P}\{V_k(Z_0\cap L)\in V_k(\varrho B_L)(1,1+h)\}
\\
&&\qquad\ge\mathbf{P}\bigl\{(X\cap L)\bigl(\mathcal{H}_{(1+\beta)\varrho
B_L}\bigr)=N(L_j^*),
Z_0\bigl((X\cap L)\,\rule{.1mm}{.26cm}\rule{.24cm}{.1mm}\, \mathcal{H}_{(1+\beta)\varrho B_L}\bigr)\in\mathcal
{P}_L\bigr\}
\\
&&\qquad=\frac{[2k(1+\beta)\varrho V_k(B_L)\gamma_{X\cap L}]^{N(L_j^*)}}{N(L_j^*)!}
\exp[-2k(1+\beta)\varrho V_k(B_L)\gamma_{X\cap L}]
\\
&&{}\qquad\quad\times  \mathbf{P}\bigl\{Z_0\bigl((X\cap L)\,\rule{.1mm}{.26cm}\rule{.24cm}{.1mm}\, \mathcal
{H}_{(1+\beta
)\varrho B_L}\bigr)\in\mathcal{P}_L
\vert(X\cap L)\bigl(\mathcal{H}_{(1+\beta)\varrho B_L}\bigr)=N(L_j^*)\bigr\}.
\end{eqnarray*}
Using a fundamental property of Poisson processes (cf. \cite{SW08}, Theorem 3.2.2(b)), the relation
$\mathbf{E}[(X\cap L)(\mathcal{H}_{(1+\beta)\varrho B_L})]=2\gamma
_{X\cap
L}k(1+\beta)\varrho V_k(B_L)$, and the definition of the set $\mathcal
{P}_L$, we obtain
\begin{eqnarray*}
&& \mathbf{P}\{V_k(Z_0\cap L)\in a(1,1+h)\}
\\
&&\qquad \ge\frac{(2\gamma_{X\cap L})^{N(L_j^*)}}{N(L_j^*)!} \exp
[-2k(1+\beta)\varrho V_k(B_L)\gamma_{X\cap L}]
\\
&& {}\qquad\quad\times  \int_{U_{N(L_j^*)}^L(L_j^*)}\cdots\int
_{U_{1}^L(L_j^*)}\int_\mathbb{R}\cdots\int_\mathbb{R}
\\
&& \qquad\quad\mathbf{1}\bigl\{|t_i-\varrho t_i^0(L_j^*)|<\varrho
\alpha(L_j^*),  \mbox{for }i=1,\ldots,N(L_j^*),\mbox{ and }
\\
&&\qquad\qquad  V_k\bigl(P\bigl(\mathbf{u}_1,\ldots,\mathbf
{u}_{N(L_j^*)};t_1,\ldots
,t_{N(L_j^*)}\bigr)\bigr)\in V_k(\varrho B_L)(1,1+h)\bigr\}
\\
&& \qquad\quad  \mathrm{d}t_1\,\cdots\, \mathrm
{d}t_{N(L_j^*)}\,
S^L_{k-1}(B_L,\mathrm{d}\mathbf{u}_1)\cdots S^L_{k-1}\bigl(B_L,\mathrm
{d}\mathbf{u}_{N(L_j^*)}\bigr),
\end{eqnarray*}
and hence
\begin{eqnarray*}
&&\mathbf{P}\{V_k(Z_0\cap L)\in a(1,1+h)\}
\\
&&\qquad\ge\frac{(2\gamma_{X\cap L})^{N(L_j^*)}}{N(L_j^*)!}\exp
[-2(1+\beta)\gamma_{X\cap L}\tau_L
a^{1/k}]\\
&&\qquad{}\quad\times c(\beta,\varphi)\varrho h(2\varrho\alpha
(L_j^*))^{N(L_j^*)-1}\prod_{i=1}^{N(L_j^*)}
S_{k-1}^L(B_L,U_i^L(L_j^*)).
\end{eqnarray*}
Since $a\ge1$, $V_k(B_L)\le c_7^k\kappa_k$ and $\gamma_{X\cap L}\ge
\gamma m(\varphi)$ [cf. (\ref{spaet2})], we finally get
\begin{eqnarray*}
&&\mathbf{P}\{V_k(Z_0\cap L)\in a(1,1+h)\}
\\
&&\qquad\ge\frac{(2a^{1/k}\gamma_{X\cap L}V_k(B_L)^{-1/k})^{N(L_j^*)}}{N(L_j^*)!}
(2\alpha(L_j^*))^{N(L_j^*)-1}c(\beta,\varphi)\eta_0^{N(L_j^*)}
\\
&&{}\qquad\quad\times  h\exp[-2(1+\beta)\gamma_{X\cap L}\tau_L
a^{1/k}]\\
&&\qquad\ge c_{10}h\exp[-2(1+\beta)\gamma_{X\cap L}\tau_L
a^{1/k}],
\end{eqnarray*}
which gives the required estimate.
\end{pf}

\begin{lemma}\label{L7}Let $0<\varepsilon<1$ and $h\in(0,1/2)$. There
are a constant $c_{11}>0$, depending only on $m(\varphi), \gamma$,
and $\varepsilon$, and a constant $c_{12}>0$, depending only on
$m(\varphi)$, such that, for $L\in G(d,k)$ and $a\ge1$,
\begin{eqnarray*}
& & \mathbf{P}\{\vartheta(Z_0 \cap L,B_L) \ge\varepsilon, V_k(Z_0
\cap L)
\in a(1,1+h)\}
\\
& &\qquad \le c_{11} h \exp[-2(1+c_{12}\varepsilon^{k+1})\gamma_{X\cap
L}\tau_L a^{1/k}].
\end{eqnarray*}
\end{lemma}

\begin{pf}The assertion is obtained by applying Proposition 7.1 of
\cite{HRS04a} in a given subspace $L\in G(d,k)$, again to a stationary
Poisson hyperplane process with intensity $\gamma_{X\cap L}$ and
spherical directional distribution $\varphi_{X\cap L}$. The slightly
different definition of the deviation measure $r_B$, as opposed to
$\vartheta$ in the present paper, is inessential for the proof. Where
a constant in \cite{HRS04a} depends on $B$, it depends now on $B_L$.
Whenever a constant in \cite{HRS04a} depends on $B$, this dependence
is via mixed volumes of $B$ with specific convex bodies, or via the
diameter of $B$, and the constant can, therefore, be estimated from the
appropriate side by positive constants for which the dependence on $B$
is only a dependence on the inradius and circumradius of $B$. Due to
the universal bounds for the inradius and circumradius of $B_L$
provided by Lemma \ref{L3}, for the constants appearing in the
application of \cite{HRS04a} to $L$, the dependence on $B_L$ is, in
fact, a dependence on $m(\varphi)$ only.
\end{pf}

\section[Proof of Theorem 2.1]{Proof of Theorem \protect\ref{T1}}\label{sec5}

Let $L^*\in G(d,k)$ with $L^*\in\operatorname{supp} \mathbf{ Q}_{d-k}$ be given.
Let $N^* \subset G(d,k)$ be a neighbourhood of $L^*$. Then
\[
\mathbf{P}\bigl\{V_k\bigl(Z_0^{(k)}\bigr)\ge a,  D\bigl(Z_0^{(k)}\bigr)\in N^*\bigr\}>0.
\]
The positivity of this probability follows from (\ref{n1}) together
with the facts that $\mathbf{ Q}_{d-k}(N^*)>0$ and that, for any $r>0$,
\[
\mathbf{P}\{ r{\mathbb B}^d\subset Z_0\}= \mathbf
{P}\{H\cap
r{\mathbb B}^d = \varnothing\enspace\forall H\in X\}>0.
\]
Let $\varepsilon>0$ and $a\ge1$. We have
\begin{eqnarray}\label{3}
& & \mathbf{P}\bigl\{\vartheta\bigl(Z_0^{(k)},B_{L^*}\bigr)\ge\varepsilon\vert
V_k\bigl(Z_0^{(k)}\bigr)\ge a,
D\bigl(Z_0^{(k)}\bigr)\in N^*\bigr\} \nonumber
\\[-8pt]\\[-8pt]
&&\qquad = \frac{\mathbf{P}\{\vartheta(Z_0^{(k)},B_{L^*})\ge
\varepsilon,
V_k(Z_0^{(k)})\ge a,
D(Z_0^{(k)})\in N^*\}}
{\mathbf{P}\{V_k(Z_0^{(k)})\ge a,  D(Z_0^{(k)})\in N^*\}}.\nonumber
\end{eqnarray}
In order to estimate this ratio, we derive an estimate from above for
the numerator and
an estimate from below for the denominator.
As in \cite{HRS04a}, we first consider the condition
$V_k(Z_0^{(k)})\in a(1,1+h)$ for $h>0$,
instead of $V_k(Z_0^{(k)})\ge a$.

For the estimate of the numerator of (\ref{3}), we use (\ref{n1}) to get
\begin{eqnarray*}
& & \mathbf{P}\bigl\{\vartheta\bigl(Z_0^{(k)},B_{L^*}\bigr)\ge\varepsilon,
V_k\bigl(Z_0^{(k)}\bigr)\in a(1,1+h),  D\bigl(Z_0^{(k)}\bigr)\in N^*\bigr\}
\\
& &\qquad = \int_{G(d,k)} \mathbf{P}\bigl\{\vartheta(Z_0\cap L,B_{L^*})\ge
\varepsilon,
V_k(Z_0\cap L)\in a(1,1+h),
\\
&&\qquad \hspace*{53pt}\hspace*{133pt}D(Z_0\cap L)\in N^*\bigr\} \mathbf{ Q}_{d-k}(\mathrm{d}L)
\\
& &\qquad = \int_{N^*} \mathbf{P}\{\vartheta(Z_0\cap L,B_{L^*})\ge
\varepsilon,
V_k(Z_0\cap L)\in a(1,1+h)\} \mathbf{ Q}_{d-k}(\mathrm{d}L).
\end{eqnarray*}

In contrast to the case of the zero cell $Z_0$ treated in \cite
{HRS04a}, we are here faced with the problem that the random polytope
$Z_0\cap L$, for variable $L$, must be compared with the fixed Blaschke
body $B_{L^*}$. This explains the necessity of restricting the
direction space $D(Z_0^{(k)})$ to a neighbourhood of $L^*$ and of
establishing the stability result Lemma \ref{L4}, which allows us the
estimate (\ref{4b}) and finally (\ref{4.6}). A similar remark
concerns the estimation of the denominator.

We choose numbers $1/2\le p < 1$ and $q>1$, depending only on $\varepsilon
$ and the number $c_{12}$ from Lemma \ref{L7} (but with $\varepsilon$
replaced by $\varepsilon/2$), such that
%
\begin{equation}\label{1b}
\frac{q}{p} < 1+\frac{c_{13}}{2}\varepsilon^{k+1}
\end{equation}
with $c_{13}:= c_{12}/2^{k+1}$. Then we choose a number $\theta>0$
satisfying the conditions
\[
\theta\le\min\biggl\{\frac{1}{8},\biggl(\frac{\varepsilon
}{2c_{9}}\biggr)^{3k}\biggr\},
\]
where $c_{9}$ is the constant from Lemma \ref{L4}, and
%
\begin{equation}\label{4b}
p \gamma_{X\cap L^*}\tau_{L^*} \le\gamma_{X\cap L}\tau_L \le q
\gamma_{X\cap L^*}\tau_{L^*} \qquad\mbox{if } \Delta(L,L^*)\le
\theta.
\end{equation}
The latter is possible by (\ref{4.3}) and Lemma \ref{L4}, since $\tau
_L= kV_k(B_L)^{1-1/k}$.

If $L\in G(d,k)$ and $\Delta(L,L^*) \le\theta$, then every convex
body $K\in\mathcal{K}_0(L)$ with $\vartheta(K,B_L) < \varepsilon/2$
satisfies $\vartheta(K,B_{L^*}) < \varepsilon$, by Lemma \ref{L5}. Now
we choose for $N^*$ the neighbourhood $ N_\theta:=N_{\theta}(L^*)$. Then
\[
L\in N_\theta\quad\mbox{and}\quad\vartheta(Z_0\cap L,B_{L^*}) \ge
\varepsilon\quad\mbox{implies}\quad\vartheta(Z_0\cap L,B_L) \ge
\varepsilon/2.
\]
This gives
\begin{eqnarray*}
& & \mathbf{P}\bigl\{\vartheta\bigl(Z_0^{(k)},B_{L^*}\bigr)\ge\varepsilon,
V_k\bigl(Z_0^{(k)}\bigr)\in a(1,1+h),  D\bigl(Z_0^{(k)}\bigr)\in N_{\theta}\bigr\}
\\
& &\qquad \le\int_{N_\theta} \mathbf{P}\{\vartheta(Z_0\cap L,B_L)\ge
\varepsilon
/2, V_k(Z_0\cap L)\in a(1,1+h)\} \mathbf{ Q}_{d-k}(\mathrm{d}L).
\end{eqnarray*}

Let $h\in(0,1/2)$. By Lemma \ref{L7} (with $\varepsilon$ replaced by
$\varepsilon/2$),
\begin{eqnarray*}
& & \mathbf{P}\{\vartheta(Z_0\cap L,B_L)\ge\varepsilon/2, V_k(Z_0\cap
L)\in
a(1,1+h)\}
\\
& &\qquad \le c_{14}h\exp[-2(1+c_{13}\varepsilon^{k+1})\gamma_{ X\cap
L}\tau_L a^{1/k}]
\end{eqnarray*}
with a constant $c_{14}$ depending only on $\varphi,\gamma,\varepsilon
$; here $c_{13}$ (defined above) depends only on $\varphi$.

By (\ref{4b}), we can conclude that
\begin{eqnarray}\label{4.6}
& & \mathbf{P}\bigl\{\vartheta\bigl(Z_0^{(k)},B_{L^*}\bigr)\ge\varepsilon,
V_k\bigl(Z_0^{(k)}\bigr)\in a(1,1+h),  D\bigl(Z_0^{(k)}\bigr)\in N_\theta\bigr\}\nonumber
\\[-8pt]\\[-8pt]
& &\qquad \le\mathbf{ Q}_{d-k}(N_\theta) c_{14}h\exp[-2(1+c_{13}\varepsilon
^{k+1})p  \gamma_{ X\cap L^*}\tau_{L^*}a^{1/k}].\nonumber
\end{eqnarray}
Now the argument in \cite{HRS04a}, pages 1164--1165 (Case 2),
leads from (\ref{4.6}) to the estimate
%
\begin{eqnarray}\label{eq24}
& & \mathbf{P}\bigl\{\vartheta\bigl(Z_0^{(k)},B_{L^*}\bigr)\ge\varepsilon,
V_k\bigl(Z_0^{(k)}\bigr)\ge a,  D\bigl(Z_0^{(k)}\bigr)\in N_\theta\bigr\} \nonumber
\\
& & \qquad \le c_{15}\mathbf{ Q}_{d-k}(N_\theta) h\exp\biggl[-2\biggl(1+\frac
{c_{13}}{2}\varepsilon^{k+1}\biggr)p  \gamma_{ X\cap L^*}\tau
_{L^*}a^{1/k}\biggr]
\\
& & {}\qquad\quad \times \exp\biggl[-\frac{c_{13}}{2}\varepsilon
^{k+1}p  \gamma_{ X\cap L^*}\tau_{L^*}a^{1/k}\biggr],\nonumber
\end{eqnarray}
where $c_{15}$ is a positive constant depending only on $\varphi,
\gamma, \varepsilon$. Here, we use that $L\mapsto\gamma_{X\cap L}$ and
$L\mapsto\tau_L$ are continuous and can be estimated from below by a
positive constant independent of $L$.

For the denominator of (\ref{3}), we obtain similarly
\begin{eqnarray*}
& & \mathbf{P}\bigl\{V_k\bigl(Z_0^{(k)}\bigr)\in a(1,1+h),  D\bigl(Z_0^{(k)}\bigr)\in
N_\theta\bigr\}
\\
& &\qquad = \int_{N_\theta} \mathbf{P}\{V_k(Z_0\cap L)\in a(1,1+h)\} \mathbf{
Q}_{d-k}(\mathrm{d}L).
\end{eqnarray*}
We define the number $\beta$, depending only on $\varphi$ and
$\varepsilon$, by
%
\begin{equation}\label{4.8}
\biggl(1+\frac{c_{13}}{2}\varepsilon^{k+1}\biggr)p=(1+\beta)q.
\end{equation}
It follows from (\ref{1b}) that $\beta>0$. By Lemma \ref{L6}, there
are constants $c_{10}$, $0<h_0<1/2$, depending only on $\varphi,
\gamma$ and $\varepsilon$, such that, for $L\in G(d,k)$, $a\ge1$ and
$0<h\le h_0$,
\[
\mathbf{P}\{V_k(Z_0\cap L) \in a(1,1+h)\} \ge c_{10} h \exp
[-2(1+\beta
)\gamma_{X\cap L}\tau_L a^{1/k}].
\]
Using (\ref{4b}) for $L\in N_\theta$, we deduce that
\begin{eqnarray*}
& & \mathbf{P}\bigl\{V_k\bigl(Z_0^{(k)}\bigr)\in a(1,1+h),  D\bigl(Z_0^{(k)}\bigr)\in
N_\theta
\bigr\} \nonumber
\\
& & \qquad\ge\mathbf{ Q}_{d-k}(N_\theta) c_{10}h\exp[-2(1+\beta)q
\gamma_{ X\cap L^*}\tau_{L^*}a^{1/k}].
\end{eqnarray*}
With $\beta$ given by (\ref{4.8}), this yields
\begin{eqnarray*}
& & \mathbf{P}\bigl\{V_k\bigl(Z_0^{(k)}\bigr)\ge a,  D\bigl(Z_0^{(k)}\bigr)\in N_\theta
\bigr\}
\\
& &\qquad \ge c_{10}\mathbf{ Q}_{d-k}(N_\theta) h\exp\biggl[-2\biggl(1+\frac
{c_{13}}{2}\varepsilon^{k+1}\biggr)p  \gamma_{ X\cap L^*}\tau
_{L^*}a^{1/k}\biggr].
\end{eqnarray*}
Here and in \eqref{eq24}, we choose the same number $h\in(0,h_0]$.
Then division gives the assertion of Theorem \ref{T1}, since $p\ge
1/2$ and we can estimate $\gamma_{ X\cap L^*}\tau_{L^*}$ from below
by a constant depending only on $\varphi$ and $\gamma$.

\section[Proof of Theorem 2.2]{Proof of Theorem \protect\ref{T2}}

The proof of Theorem \ref{T2} is based on \eqref{1a}, which is
applied with different functions $f$, and on the relation
\[
\mathbf{E}V_k \bigl(Z^{(k)}\bigr)=\frac{d^{(k)}_k}{\gamma^{(k)}}=\frac
{V_{d-k}(\Pi
_X)}{\left({d\atop k}\right)V_d(\Pi_X)}
=:c_{16},
\]
which follows from \cite{SW08}, equation (10.3) and Theorem 10.3.3,
with $c_{16}$ depending only on $\varphi$ and $\gamma$.

We use definitions and results from the preceding proof of Theorem \ref
{T1}. In particular, $\beta$ is defined by (\ref{4.8}). Then there
are positive constants $c_{17}, \theta_1$ and $h_1 < 1/2$, depending
only on $\varphi,\gamma$ and $\varepsilon$, such that, for $a\ge1$,
$0<\theta\le\theta_1$ and $0<h\le h_1$,
\begin{eqnarray*}
&&\mathbf{P}\bigl\{V_k\bigl(Z_0^{(k)}\bigr)\in a(1,1+h), D\bigl(Z_0^{(k)}\bigr)\in
N_{\theta}\bigr\}
\\
&&\qquad\ge\mathbf{ Q}_{d-k}(N_{\theta})c_{17}h\exp\biggl[-2\biggl(1+\frac
{\beta}{2}\biggr)
q \gamma_{X\cap L^*}\tau_{L^*}a^{1/k}\biggr].
\end{eqnarray*}
For a polytope $K\subset\mathbb{R}^d$, we now define
\[
f(K):=\mathbf{1}\{V_k(K)\in a(1,1+h), D(K)\in N_{\theta}\}V_k(K)^{-1},
\]
if $K$ is $k$-dimensional, and $f(K):=0$ otherwise. Clearly, $f$ is
translation invariant, and for $a\ge1$ and $0<h\le h_1$, \eqref{1a} gives
\begin{eqnarray*}
&&\mathbf{P}\bigl\{V_k\bigl(Z^{(k)}\bigr)\in a(1,1+h), D\bigl(Z^{(k)}\bigr)\in N_{\theta}\bigr\}
\\
&&\qquad=\mathbf{E}V_k\bigl(Z^{(k)}\bigr)\mathbf{E}\bigl[\mathbf{1}\bigl\{
V_k\bigl(Z^{(k)}_0\bigr)\in a(1,1+h),
D\bigl(Z^{(k)}_0\bigr)\in N_{\theta}\bigr\}V_k\bigl(Z^{(k)}_0\bigr)^{-1}\bigr]
\\
&&\qquad\ge c_{16}\frac{1}{1+h_1}\frac{1}{a}\mathbf{P}\bigl\{
V_k\bigl(Z^{(k)}_0\bigr)\in
a(1,1+h), D\bigl(Z^{(k)}_0\bigr)\in N_{\theta}\bigr\}
\\
&&\qquad\ge c_{18}\mathbf{ Q}_{d-k}(N_{\theta})\frac{1}{a}h
\exp\biggl[-2\biggl(1+\frac{\beta}{2}\biggr)q \gamma_{X\cap
L^*}\tau_{L^*}a^{1/k}\biggr]
\\
&&\qquad\ge c_{19}\mathbf{ Q}_{d-k}(N_{\theta})h
\exp[-2(1+\beta)q \gamma_{X\cap L^*}\tau
_{L^*}a^{1/k}],
\end{eqnarray*}
since $\gamma_{X\cap L^*}\tau_{L^*}\ge c_{20}>0$. Here, $c_{18}$ and
$c_{19}$ depend only on
$\varphi,\gamma,\varepsilon$, and $c_{20}$ depends only on $\varphi
,\gamma$.
In particular, recalling the definition of $\beta$ from (\ref{4.8}),
\begin{eqnarray}\label{eqA}
&&\mathbf{P}\bigl\{V_k\bigl(Z^{(k)}\bigr)\ge a,  D\bigl(Z^{(k)}\bigr)\in N_{\theta
}\bigr\}\nonumber
\\[-8pt]\\[-8pt]
&&\qquad\ge c_{19}\mathbf{ Q}_{d-k}(N_{\theta})h_1
\exp\biggl[-2\biggl(1+\frac{c_{13}}{2}\varepsilon^{k+1}\biggr)p
\gamma_{X\cap L^*}\tau_{L^*}a^{1/k}\biggr].\nonumber
\end{eqnarray}

For the upper bound, we put
\[
f(K):=\mathbf{1}\{\vartheta(K,B_{L^*})\ge\varepsilon,  V_k(K)\ge a,
D(K)\in N_{\theta}\}V_k(K)^{-1},
\]
if $K$ is a $k$-dimensional polytope, and $f(K):=0$ otherwise, where
$0<\theta\le\theta_1$, with $\theta_1$ sufficiently small, and
$a\ge1$. Using again \eqref{1a}, we obtain
\begin{eqnarray}\label{eqB}
&&\mathbf{P}\bigl\{\vartheta\bigl(Z^{(k)},B_{L^*}\bigr)\ge\varepsilon,
V_k\bigl(Z^{(k)}\bigr)\ge a,  D\bigl(Z^{(k)}\bigr)\in N_{\theta}\bigr\}\nonumber
\\
&&\qquad=c_{16}\mathbf{E}\bigl[\mathbf{1}\bigl\{\vartheta
\bigl(Z^{(k)}_0,B_{L^*}\bigr)\ge
\varepsilon,  V_k\bigl(Z^{(k)}_0\bigr)\ge a,   D\bigl(Z^{(k)}_0\bigr)\in N_{\theta}
\bigr\}V_k\bigl(Z_0^{(k)}\bigr)^{-1}\bigr]\nonumber
\\[-8pt]\\[-8pt]
&&\qquad\le c_{21}\mathbf{ Q}_{d-k}(N_{\theta})h_1
\exp\biggl[-2\biggl(1+\frac{c_{13}}{2}\varepsilon^{k+1}\biggr)p
\gamma_{X\cap L^*}\tau_{L^*}a^{1/k}\biggr]\nonumber
\\
&&\qquad\quad{}\times\exp\biggl[-\frac{c_{13}}{2}\varepsilon^{k+1}p
\gamma_{X\cap L^*}\tau_{L^*}a^{1/k}\biggr],\nonumber
\end{eqnarray}
where \eqref{eq24} was used in the last estimate and $c_{21}$ depends
only on $\varphi,\gamma,\varepsilon$.

From \eqref{eqA} and \eqref{eqB}, we conclude that
\begin{eqnarray*}
&&\mathbf{P}\bigl\{\vartheta\bigl(Z^{(k)},B_{L^*}\bigr)\ge\varepsilon\vert
V_k\bigl(Z^{(k)}\bigr)\ge a,  D\bigl(Z^{(k)}\bigr)\in N_{\theta}\bigr\}
\\
&&\qquad\le c_{22}\exp\biggl[-\frac{c_{13}}{2}\varepsilon^{k+1}p \gamma
_{X\cap L^*}\tau_{L^*}a^{1/k}\biggr]
\\
&&\qquad\le c_{22}\exp[-c_{23}\varepsilon^{k+1}a^{1/k}],
\end{eqnarray*}
where $c_{22}$ depends only on $\varphi,\gamma,\varepsilon$ and
$c_{23}$ depends only on $\varphi$ and $\gamma$.

\section{Limit shapes}\label{sec7}

Similarly as in \cite{HS07a}, Section 4, but with an additional limit
procedure referring to direction spaces, we can establish the existence
of limit shapes.

For a convex body $K\subset{\mathbb R}^d$, we denote by $s_{\sf H}(K)$
the equivalence class of all convex bodies homothetic to $K$; this is
the (homothetic) \textit{shape} of $K$. Let $\mathcal{S}_{\sf H}$ denote the
space of all shapes, equipped with the quotient topology.

Let the assumptions of Theorem \ref{T2} be satisfied; in particular,
$L^*\in G(d,k)$ is contained in the support of the measure $\mathbf{ Q}_{d-k}$.

The \textit{conditional law of the shape of} $Z_0^{(k)}$, given the lower
bound $a$ for its $k$-volume and the upper bound $\theta$ for the
distance of its direction space from $L^*$, is defined by
\[
\mu_{a,\theta}(A):= \mathbf{P}\bigl\{s_{\sf H}\bigl(Z_0^{(k)}\bigr)\in A\vert
V_k\bigl(Z_0^{(k)}\bigr) \ge a, \Delta\bigl(D\bigl(Z_0^{(k)}\bigr),L^*\bigr) <\theta\bigr\}
\]
for $A\in\mathcal{B}(\mathcal{S}_{\sf H})$.

\begin{theorem} \label{T7}
The shape $s_{\sf H}(B_{L^*})$ is the limit shape of the weighted
typical cell $Z_0^{(k)}$ with respect to $V_k$ and $\Delta(D(\cdot
),L^*)$, in the sense that
\[
\lim_{a\to\infty\atop\theta\to0} \mu_{a,\theta} = \delta
_{s_{\sf H}(B_{L^*})}\qquad\mbox{weakly},
\]
where $\delta_{s_{\sf H}(B_{L^*})}$ denotes the Dirac measure
concentrated at $s_{\sf H}(B_{L^*})$.
\end{theorem}

\begin{pf}Let $\mathcal{C}\subset\mathcal{S}_{\sf H}$ be closed. It
suffices to show that
%
\begin{equation}\label{7.1}
\mathop{\limsup_{a\to\infty}}_{\theta\to0} \mu_{a,\theta}(\mathcal{C})
\le\delta_{s_{\sf H}(B_{L^*})}(\mathcal{C}).
\end{equation}
We assume that $s_{\sf H}(B_{L^*}) \notin\mathcal{C}$ and that $\mathcal{C}$
contains the shape of at least one $k$-dimensional body, since
otherwise (\ref{7.1}) holds trivially. For $K\in\mathcal{K}$ with $\dim
K=k$, we put $f(K):=\vartheta(K,B_{L^*})+\Delta(D(K),L^*)$. Let
\begin{eqnarray*}
\mathcal{K}^* &:=& \{ K \in\mathcal{K}\dvtx \dim K=k, s_{\sf H}(K) \in\mathcal{
C},  B_{D(K)}\subset K\},
\\
\alpha &:=& \inf_{K\in\mathcal{K}^*}f(K),
\end{eqnarray*}
and choose $c>\alpha$. There exists $R>0$ such that every $K\in\mathcal{
K}^*$ with $f(K)\le c$ has a homothetic copy that is contained in
$R{\mathbb B}^d$. Hence, if we put
\[
\mathcal{K}_c^*:=\{K\in\mathcal{K}^*\dvtx f(K)\le c,  K\subset R{\mathbb B}^d\},
\]
then $\alpha= \inf_{K\in\mathcal{K}_c^*}f(K)$. The function $f$ is
continuous and the set $\mathcal{K}_c^*$ is compact (note that the
condition $B_{D(K)}\subset K$ in the definition of $\mathcal{K}^*$ ensures
that limits of bodies in $\mathcal{K}^*$ still have dimension $k$).
Therefore, the infimum $\alpha$ is attained, say at $K_0$. If $\alpha
=0$, then $K_0$ is homothetic to $B_{L^*}$, hence $s_{\sf H}(B_{L^*}) =
s_{\sf H}(K_0)\in\mathcal{C}$, a contradiction. It follows that $\alpha>0$.

Put $\varepsilon:=\alpha/2$. To this $\varepsilon$, we can choose
constants $c_1,c_2,c_3$ according to Theorem~2.1, such that
\begin{eqnarray*}
& & \mathbf{P}\bigl\{ \vartheta\bigl(Z_0^{(k)},B_{L^*}\bigr) \ge\varepsilon\vert
V_k\bigl(Z_0^{(k)}\bigr) \ge a,  D\bigl(Z_0^{(k)}\bigr) \in N_{\theta}(L^*)\bigr\}
\\
& &\qquad \le c_2 \exp[-c_3 \varepsilon^{k+1}a^{1/k}]
\end{eqnarray*}
for $a\ge1$ and $0<\theta\le c_1$.

Every $k$-dimensional convex body $K\in s_{\sf H}^{-1}(\mathcal{C})$ with
$\Delta(D(K),L^*)\le\alpha/2$ satisfies $\vartheta(K,B_{L^*}) \ge
\varepsilon$. Hence, for $0<\theta\le\min\{c_1,\alpha/2\}$ we have
\begin{eqnarray*}
\mu_{a,\theta}(\mathcal{C}) &=& \mathbf{P}\bigl\{ s_{\sf
H}\bigl(Z_0^{(k)}\bigr)\in
\mathcal{C}\vert V_k\bigl(Z_0^{(k)}\bigr)\ge a,  D\bigl(Z_0^{(k)}\bigr)\in N_\theta
(L^*)\bigr\}
\\
&\le& \mathbf{P}\bigl\{\vartheta\bigl(Z_0^{(k)}, B_{L^*}\bigr)\ge\varepsilon
\vert
V_k\bigl(Z_0^{(k)}\bigr)\ge a,  D\bigl(Z_0^{(k)}\bigr)\in N_\theta(L^*)\bigr\}
\\
&\le& c_2 \exp[-c_3 \varepsilon^{k+1}a^{1/k}].
\end{eqnarray*}
For $a\to\infty$ this tends to zero, hence (\ref{7.1}) follows.
\end{pf}

Theorem 2.2 yields a completely analogous result for the typical cell.

\section*{Acknowledgment} We thank the referee for his/her very careful
reading of the manuscript and for several valuable suggestions for improvements.

%

\printaddresses


\begin{thebibliography}{17}

\bibitem{BL07}
%
\begin{barticle}[mr]
\bauthor{\bsnm{Baumstark},~\bfnm{Volker}\binits{V.}} \AND
\bauthor{\bsnm{Last},~\bfnm{G{\"u}nter}\binits{G.}}
(\byear{2007}).
\btitle{Some distributional results for {P}oisson--{V}oronoi tessellations}.
\bjournal{Adv. in Appl. Probab.}
\bvolume{39}
\bpages{16--40}.
\bid{doi={10.1239/aap/1175266467}, mr={2307869}}
\end{barticle}
%
\endbibitem

\bibitem{DK70}
%
\begin{barticle}[mr]
\bauthor{\bsnm{Davis},~\bfnm{Chandler}\binits{C.}} \AND
\bauthor{\bsnm{Kahan},~\bfnm{W.~M.}\binits{W.~M.}}
(\byear{1970}).
\btitle{The rotation of eigenvectors by a perturbation. {III}}.
\bjournal{SIAM J. Numer. Anal.}
\bvolume{7}
\bpages{1--46}.
\bid{mr={0264450}}
\end{barticle}
%
\endbibitem

\bibitem{GKW08}
%
\begin{bmisc}[auto:SpringerTagBib|2009-01-14|16:51:27]
\bauthor{\bsnm{Goodey},~\bfnm{P.}\binits{P.}},
\bauthor{\bsnm{Kiderlen},~\bfnm{M.}\binits{M.}}
\AND
\bauthor{\bsnm{Weil},~\bfnm{W.}\binits{W.}}
(\byear{2010}).
\bhowpublished{Spherical projections and liftings in geometric
tomography. \textit{Adv. Geom.} To appear}.
\end{bmisc}
%
\endbibitem

\bibitem{HRS04a}
%
\begin{barticle}[mr]
\bauthor{\bsnm{Hug},~\bfnm{Daniel}\binits{D.}},
\bauthor{\bsnm{Reitzner},~\bfnm{Matthias}\binits{M.}} \AND
\bauthor{\bsnm{Schneider},~\bfnm{Rolf}\binits{R.}}
(\byear{2004}).
\btitle{The limit shape of the zero cell in a stationary {P}oisson hyperplane
tessellation}.
\bjournal{Ann. Probab.}
\bvolume{32}
\bpages{1140--1167}.
\bid{doi={10.1214/aop/1079021474}, mr={2044676}}
\end{barticle}
%
\endbibitem

\bibitem{HRS04b}
%
\begin{barticle}[mr]
\bauthor{\bsnm{Hug},~\bfnm{Daniel}\binits{D.}},
\bauthor{\bsnm{Reitzner},~\bfnm{Matthias}\binits{M.}} \AND
\bauthor{\bsnm{Schneider},~\bfnm{Rolf}\binits{R.}}
(\byear{2004}).
\btitle{Large {P}oisson--{V}oronoi cells and {C}rofton cells}.
\bjournal{Adv. in Appl. Probab.}
\bvolume{36}
\bpages{667--690}.
\bid{doi={10.1239/aap/1093962228}, mr={2079908}}
\end{barticle}
%
\endbibitem

\bibitem{HS02}
%
\begin{barticle}[mr]
\bauthor{\bsnm{Hug},~\bfnm{Daniel}\binits{D.}} \AND
\bauthor{\bsnm{Schneider},~\bfnm{Rolf}\binits{R.}}
(\byear{2002}).
\btitle{Stability results involving surface area measures of convex bodies}.
\bjournal{Rend. Circ. Mat. Palermo (2) Suppl.}
\bvolume{70}
\bpages{21--51}.
\bid{mr={1962583}}
\end{barticle}
%
\endbibitem

\bibitem{HS04}
%
\begin{barticle}[mr]
\bauthor{\bsnm{Hug},~\bfnm{Daniel}\binits{D.}} \AND
\bauthor{\bsnm{Schneider},~\bfnm{Rolf}\binits{R.}}
(\byear{2004}).
\btitle{Large cells in {P}oisson--{D}elaunay tessellations}.
\bjournal{Discrete Comput. Geom.}
\bvolume{31}
\bpages{503--514}.
\bid{doi={10.1007/s00454-003-0818-3}, mr={2053496}}
\end{barticle}
%
\endbibitem

\bibitem{HS05}
%
\begin{barticle}[mr]
\bauthor{\bsnm{Hug},~\bfnm{Daniel}\binits{D.}} \AND
\bauthor{\bsnm{Schneider},~\bfnm{Rolf}\binits{R.}}
(\byear{2005}).
\btitle{Large typical cells in {P}oisson--{D}elaunay mosaics}.
\bjournal{Rev. Roumaine Math. Pures Appl.}
\bvolume{50}
\bpages{657--670}.
\bid{mr={2204143}}
\end{barticle}
%
\endbibitem

\bibitem{HS07a}
%
\begin{barticle}[mr]
\bauthor{\bsnm{Hug},~\bfnm{Daniel}\binits{D.}} \AND
\bauthor{\bsnm{Schneider},~\bfnm{Rolf}\binits{R.}}
(\byear{2007}).
\btitle{Asymptotic shapes of large cells in random tessellations}.
\bjournal{Geom. Funct. Anal.}
\bvolume{17}
\bpages{156--191}.
\bid{doi={10.1007/s00039-007-0592-0}, mr={2306655}}
\end{barticle}
%
\endbibitem

\bibitem{HS07b}
%
\begin{barticle}[mr]
\bauthor{\bsnm{Hug},~\bfnm{Daniel}\binits{D.}} \AND
\bauthor{\bsnm{Schneider},~\bfnm{Rolf}\binits{R.}}
(\byear{2007}).
\btitle{Typical cells in {P}oisson hyperplane tessellations}.
\bjournal{Discrete Comput. Geom.}
\bvolume{38}
\bpages{305--319}.
\bid{doi={10.1007/s00454-007-1340-9}, mr={2343310}}
\end{barticle}
%
\endbibitem

\bibitem{Kov97}
%
\begin{barticle}[mr]
\bauthor{\bsnm{Kovalenko},~\bfnm{I.~N.}\binits{I.~N.}}
(\byear{1997}).
\btitle{A proof of a conjecture of {D}avid {K}endall on the shape of random
polygons of large area}.
\bjournal{Kibernet. Sistem. Anal.}
\bvolume{4}
\bpages{3--10, 187}.
\bid{doi={10.1007/BF02733102}, mr={1609157}}
\end{barticle}
%
\endbibitem

\bibitem{Kov99}
%
\begin{barticle}[mr]
\bauthor{\bsnm{Kovalenko},~\bfnm{Igor~N.}\binits{I.~N.}}
(\byear{1999}).
\btitle{A simplified proof of a conjecture of {D}. {G}. {K}endall concerning
shapes of random polygons}.
\bjournal{J. Appl. Math. Stochastic Anal.}
\bvolume{12}
\bpages{301--310}.
\bid{mr={1736071}}
\end{barticle}
%
\endbibitem

\bibitem{PW94}
%
\begin{barticle}[mr]
\bauthor{\bsnm{Paige},~\bfnm{C.~C.}\binits{C.~C.}} \AND
\bauthor{\bsnm{Wei},~\bfnm{M.}\binits{M.}}
(\byear{1994}).
\btitle{History and generality of the {${\rm CS}$} decomposition}.
\bjournal{Linear Algebra Appl.}
\bvolume{208/209}
\bpages{303--326}.
\bid{doi={10.1016/0024-3795(94)90446-4}, mr={1287355}}
\end{barticle}
%
\endbibitem

\bibitem{Sch93}
%
\begin{bbook}[mr]
\bauthor{\bsnm{Schneider},~\bfnm{Rolf}\binits{R.}}
(\byear{1993}).
\btitle{Convex Bodies: The {B}runn--{M}inkowski Theory}.
\bseries{Encyclopedia of Mathematics and Its Applications}
\bvolume{44}.
\bpublisher{Cambridge Univ. Press}, \baddress{Cambridge}.
\bid{mr={1216521}}
\end{bbook}
%
\endbibitem

\bibitem{Sch09}
%
\begin{barticle}[vtex]
\bauthor{\bsnm{Schneider},~\bfnm{R.}\binits{R.}}
(\byear{2009}).
\btitle{Weighted faces of Poisson hyperplane tessellations}.
\bjournal{Adv. in Appl. Probab.}
\bvolume{41}
\bpages{682--694}.
\end{barticle}
%
\endbibitem

\bibitem{SW08}
%
\begin{bbook}[mr]
\bauthor{\bsnm{Schneider},~\bfnm{Rolf}\binits{R.}} \AND
\bauthor{\bsnm{Weil},~\bfnm{Wolfgang}\binits{W.}}
(\byear{2008}).
\btitle{Stochastic and Integral Geometry}.
\bpublisher{Springer}, \baddress{Berlin}.
\bid{mr={2455326}}
\end{bbook}
%
\endbibitem

\bibitem{SKM95}
%
\begin{bbook}[mr]
\bauthor{\bsnm{Stoyan},~\bfnm{D.}\binits{D.}},
\bauthor{\bsnm{Kendall},~\bfnm{W.~S.}\binits{W.~S.}} \AND
\bauthor{\bsnm{Mecke},~\bfnm{J.}\binits{J.}}
(\byear{1995}).
\btitle{Stochastic Geometry and Its Applications}.
\bpublisher{Wiley}, \baddress{Chichester}.
\bid{mr={895588}}
\end{bbook}
%
\endbibitem

\end{thebibliography}
\end{document}